%% file: main.tex
\documentclass[%
 onecolumn, 12pt, tightenlines
]{article}

\usepackage{amsmath, amssymb}
\usepackage[top=2cm, left=2.5cm, right=2.5cm]{geometry}
\usepackage{graphicx}
\usepackage{dcolumn}
\usepackage{bbm}
\usepackage{hyperref}
\usepackage{csquotes}
\usepackage{float}
\usepackage{caption}
\usepackage{subcaption}
\usepackage{amsthm}
\usepackage{tabularx}
\usepackage{color}
\usepackage{appendix}
\usepackage{authblk}

\newtheorem{defn}{Definition}[section]

\newtheorem{conjecture}{Conjecture}[section]

\newcommand{\E}{\mathbb{E}}
\newcommand{\prob}{\mathbb{P}}
\newcommand{\Reals}{\mathbb{R}}

\newcommand{\niso}{N_{\textrm{iso}}}
\newcommand{\nucg}{N_{\textrm{ucg}}}
\newcommand{\piso}{P_{\textrm{iso}}}
\newcommand{\pucg}{P_{\textrm{ucg}}}
\newcommand{\pfc}{P_{\textrm{dis}}}
\newcommand{\pisoucg}{P_{\textrm{iso}\cup\textrm{ucg}}}

\newcommand{\nsims}{5\,000}
\newcommand{\eqn}{\textrm{Eqn.}}
\newcommand{\eqns}{\textrm{Eqns.}}

\title{Connectivity in One-Dimensional Soft Random Geometric Graphs}
\author[ ]{Michael Wilsher \hspace{-1ex}}
\author[ ]{Carl P. Dettmann \hspace{-1ex}}
\author[ ]{Ayalvadi Ganesh \hspace{-1ex}}
\affil[ ]{\hspace{-1ex} School of Mathematics, University of Bristol, Woodland Road, Bristol, BS8 1UG, UK}
\date{\today}                     
\begin{document}

\maketitle




\input{Section0_Abstract.tex}


\input{Section1_Intro.tex}
\input{Section2_Model.tex}
\input{Section3_Analysis.tex}
\input{Section4_Iso.tex}

\input{Section6_UCG_New.tex}
\input{Section8_Conclusion.tex}
\input{Section9_Ack.tex}
\input{Section_Appendix.tex}

\bibliographystyle{ieeetr}
\bibliography{references_2.bib}

\end{document}

%% file: Section0_Abstract.tex
\begin{abstract}
In this paper, we study the connectivity of a one-dimensional soft random geometric graph (RGG). The graph is generated by placing points at random on a bounded line segment and connecting pairs of points with a probability that depends on the distance between them. We derive bounds on the probability that the graph is fully connected by analysing key modes of disconnection. In particular, analytic expressions are given for the mean and variance of the number of isolated nodes, and a sharp threshold established for their occurrence. Bounds are also derived for uncrossed gaps, and it is shown analytically that uncrossed gaps have negligible probability in the scaling at which isolated nodes appear. This is in stark contrast to the hard RGG in which uncrossed gaps are the most important factor when considering network connectivity. 
\end{abstract}

%% file: Section1_Intro.tex
\section{Introduction}

The original random geometric graph (RGG) model, also known as the unit disk or Gilbert graph or Boolean model~\cite{gilbert1961random}, was proposed as an extension of the Erd\H{o}s-R\'{e}nyi (ER) random graph~\cite{erdHos1960evolution} in which the spatial locations of nodes are taken into consideration. It is generated by distributing points (or nodes) randomly on some set, typically $\Reals^d$ or a bounded subset of it, and connecting node pairs whose distance is smaller than some threshold. Instead, if the edge between a pair of nodes is present with a probability $H(r)$ that depends on the distance $r$ between them, independent of all other edges, the resulting model has been variously termed a random connection model~\cite{meester1996continuum, mao2013connectivity}, a soft RGG~\cite{penrose2016connectivity, dettmann2016random}, or a Waxman graph~\cite{waxman1988routing}. We will call them \enquote{soft} RGGs, and term $H(\cdot)$ the connection function. The \enquote{hard} RGG is the Gilbert model, where $H(r)=1$ if $r\leq r_c$, and $H(r)=0$ otherwise; $r_c>0$ is a parameter of the model. The terms \textit{points} and \textit{nodes} will be used interchangeably henceforth.

Hard RGG models have been widely applied, e.g., to disease spread \cite{eubank2004modelling}, climate dynamics \cite{donges2009complex}, infrastructure networks \cite{robson2015resilience}, and neuronal networks \cite{nicosia2013phase}; see the survey in \cite{barthelemy2011spatial} for more examples. An extensive and rigorous mathematical study of their properties is presented in~\cite{penrose2003random}. 

The one-dimensional setting is motivated by vehicular ad-hoc networks (VANETs), which are expected to be essential for autonomous vehicles; these will be fitted with on-board radios to enable the exchange of location and velocity data, manoeuvring intentions, and safety critical data such as crash warnings. The road is modelled as a line, with nodes representing vehicles and an edge between two nodes indicating that two vehicles can communicate directly with each other. A key question is connectivity: When is every vehicle in the network (defined as a stretch of road) able to communicate with every other vehicle via a single- or multi-hop path? The hard RGG case has been studied in \cite{devroye1981laws, han2007very, ajeer2011network, knight2016counting}.

It is also of interest to identify the type of event that leads to disconnection. In two or more dimensions, the most likely cause for full connectivity to fail, in both hard and soft RGGs, is the existence of isolated nodes~\cite{penrose2003random, penrose2016connectivity}. In 1-D hard RGGs, it is the partition of the network into two or more connected clusters with no edges between them. Indeed, note that if a node $v$ on the line is isolated, then there are no edges between nodes to the left of $v$ and those to its right in the hard RGG. Such a partition can occur even without isolated nodes being present. It is not obvious which of these mechanisms dominates in soft RGGs, and has not been previously resolved. This is the main topic of this paper.

The rest of this paper is laid out as follows. The model is introduced in Section \ref{sec:model}. A precise definition of isolated nodes and uncrossed gaps is given in Section \ref{sec:analysis}, along with simulation results suggesting that these are the primary causes of disconnection. This motivates the rigorous analysis of node isolation presented in Section \ref{sec:iso}, and of uncrossed gaps in Section \ref{sec:ucg}. Sections \ref{sec:iso} and \ref{sec:ucg} also give the most significant result in this paper, namely that isolated nodes are a more important factor than uncrossed gaps when analysing network connectivity, even with a rapidly decaying connection function. Section \ref{sec:conclusions} concludes the paper with a discussion of future directions.

%% file: Section2_Model.tex
\section{Model}
\label{sec:model}

Node locations in RGGs are typically modelled by point processes~\cite{haenggi2012stochastic}, and most commonly by a Poisson process, which we define before describing the model. All sets and functions referred to in this paper are assumed to be Borel measurable, even if not explicitly stated.

\begin{defn}[Poisson Point Process]
\quad\newline
Let $\lambda: \Reals^d \to \Reals_+$ be a function whose integrals on bounded subsets of $\Reals^d$ are finite. A Poisson point process (PPP) $\Phi$ on $\Reals^d$ with intensity $\lambda(\cdot)$ is a random set of points such that $\Phi(B)$, the number of points in $B\subset \Reals^d$, has a Poisson distribution with mean $\Lambda(B)=\int_B \lambda(x)dx$, while the numbers of points in disjoint sets are mutually independent random variables.
If $\lambda(x)$ is identically equal to a constant $\lambda$, we say that the PPP is homogenous with intensity $\lambda$.
\end{defn}

\noindent \emph{Our Model} The nodes of the RGG are the points of a homogenous PPP $\Phi$ of intensity $\lambda=1$ on an interval $[0,L]$, where $L\in (0,\infty)$. The edge between points $x$ and $y$ in $\Phi$  is present with probability $H(|x-y|)$, independent of the point configuration and of all other edges. Here, $|x-y|$ denotes the Euclidean distance between $x$ and $y$, and $H:\Reals_+ \to [0,1]$ is the connection function.

The model is parametrised by the scalar $L$ and the function $H(\cdot)$. We impose mild restrictions on $H$ in Sections \ref{sec:iso} and \ref{sec:ucg}, but do not assume a specific functional form. There is no loss of generality in assuming that $\lambda=1$ as this simply defines the unit of length. Consequently, $\lambda L$ is a dimensionless quantity rather than having units of length. Another common scaling in the literature is to set $L=1$, while making $\lambda$ a free parameter. These are equivalent up to a suitable rescaling of $H(\cdot)$ \cite{mao2011towards}.

\begin{figure}
    \centering
    \includegraphics[width=0.45\textwidth]{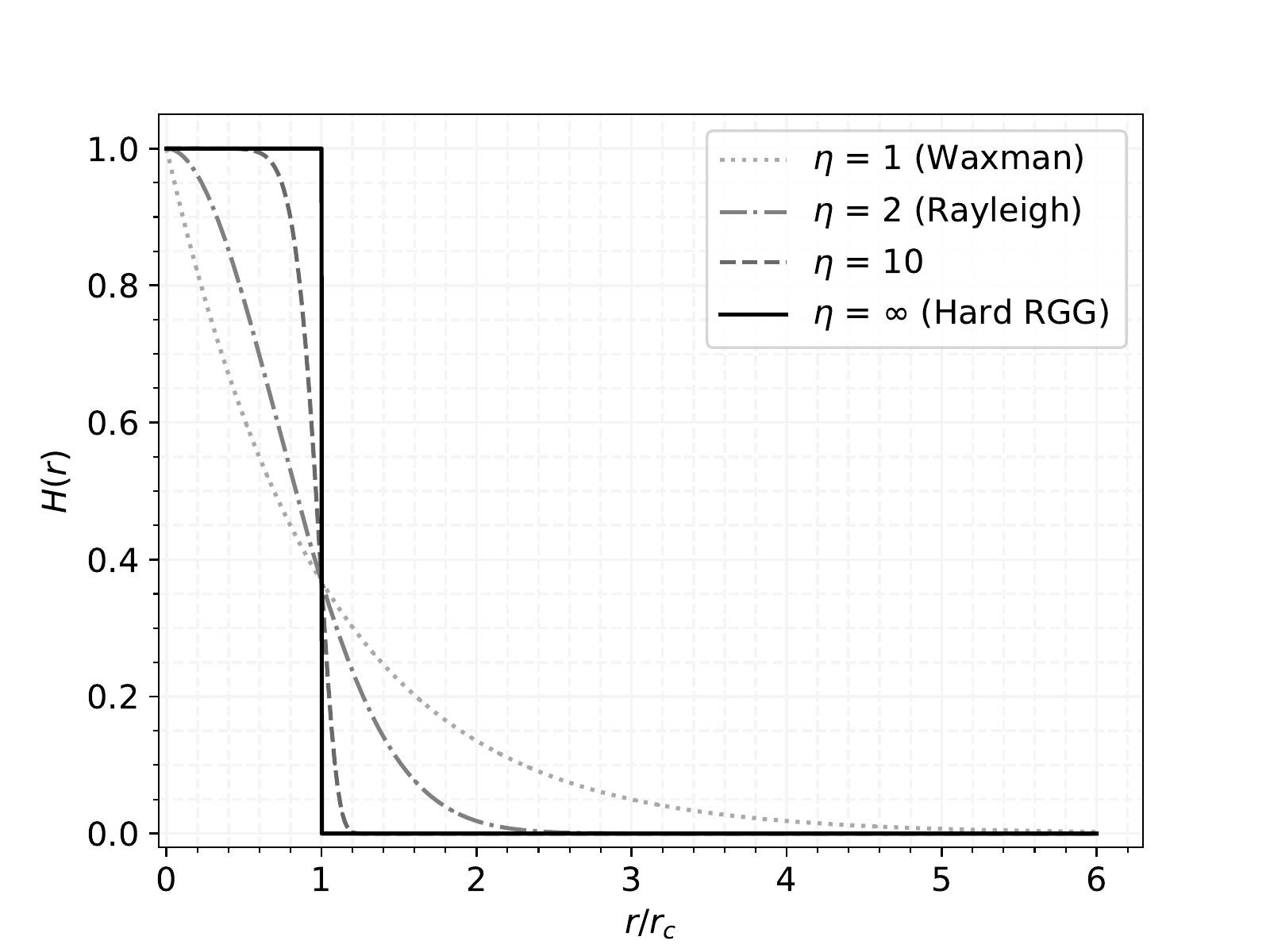}
    \caption{The probability that two nodes, a distance $r$ apart, are connected for different values of $\eta$ as defined in $\eqn$ (\ref{eqn:GeneralRayleighConn}) with $\beta=  r_c = 1$.}
    \label{fig:conn_funs}
\end{figure}

The choice of a Poisson process is made primarily for analytical tractability, but is realistic for traffic in its free flow state~\cite{roess2004traffic}, and hence for VANETs. The most widely used connection functions in the wireless communications literature are of the form 
\begin{equation}
\label{eqn:GeneralRayleighConn}
    H(r) = \beta \textrm{exp}\left( - (r/r_c)^{\eta} \right),
\end{equation}
where $r_c>0$ specifies the link range, while $\eta>0$ is related to the path loss exponent; $\beta \in (0,1]$ specifies the edge probability for short-range connections and is usually taken to be 1. The most common examples are the Waxman and Rayleigh connection functions, which correspond to $\beta = 1$, and $\eta=1$ and $\eta=2$ respectively:
\begin{equation}
\label{conn_fun_wax_ray}
    H_{\textrm{Wax}}(r) = e^{-(r/r_c)}, \quad
    H_{\textrm{Ray}}(r) = e^{-(r/r_c)^2}.
\end{equation}
The hard RGG is recovered in the limit of $\eta$ tending to infinity. Figure \ref{fig:conn_funs} depicts edge probabilites as a function of distance, for different connection functions from this class.
More general connection functions, many of which have the same general shape, may be found in \cite{Dettmann2016}.

Our goal is to derive expressions for the probability that a soft RGG generated by our model is fully connected. This appears intractable, so we restrict attention to two specific modes of disconnection, namely isolated nodes and uncrossed gaps. These are defined in the next section, and shown to account for most disconnections in extensive simulations. But even for these events, exact expressions for their probabilities cannot be obtained in closed form. Hence, we study an asymptotic regime where $L$ tends to infinity while the connection function is rescaled as $H^L(\cdot)=H(\cdot/R_L)$; we seek to identify a scaling regime $R_L$, tending to infinity at a 
specified rate with $L$, which is critical for the emergence of isolated nodes or uncrossed gaps. 

It is known in two or more dimensions that, for natural analogues of our model, the existence of isolated nodes exhibits a sharp change at a scaling of $R_L=\ln L$.  That is, for $R_L$ growing significantly faster than $\ln L$ there are longer range connections and no isolated nodes with high probability, whilst for $R_L$ growing significantly slower than $\ln L$ there are infinitely many isolated nodes in the limit. We show in Section~\ref{sec:iso} that a similar change at a $\ln L$ scaling is also observed in 1-D if the connection function is integrable. In Section~\ref{sec:ucg}, we show, under the additional condition that the connection function is monotonically decreasing and has unbounded support, that uncrossed gaps have a vanishingly small probability of occurring in this scaling regime. This is true even for rapidly decaying connection functions, and is in stark contrast to the hard RGG model, where uncrossed gaps are at least as likely as isolated nodes.

%% file: Section3_Analysis.tex
\section{Simulations}
\label{sec:analysis}

A soft RGG in 1-D may fail to be connected in many ways. We conjecture that the two main obstructions to full connectivity are the presence of isolated nodes or uncrossed gaps, defined below. In the following, $V$ denotes a realisation of $\Phi$, i.e., a point configuration comprising the node set of the graph.
\begin{defn}[Isolated Node]
    \quad \newline
    A node $x \in V$ is said to be isolated if there is no edge between the point $x$ and any point $y \in V\setminus \lbrace x\rbrace$.
\end{defn}




\begin{defn}[Uncrossed Gap]
    \quad \newline
    An uncrossed gap is said to occur at $x\in V$ if $V\cap (x,L]$ is non-empty and there are no edges between $V\cap [0,x]$ and $V\cap (x,L]$.
\end{defn}

\begin{figure}
    \centering
    \begin{subfigure}{0.45\textwidth}
    \includegraphics[width=\textwidth]{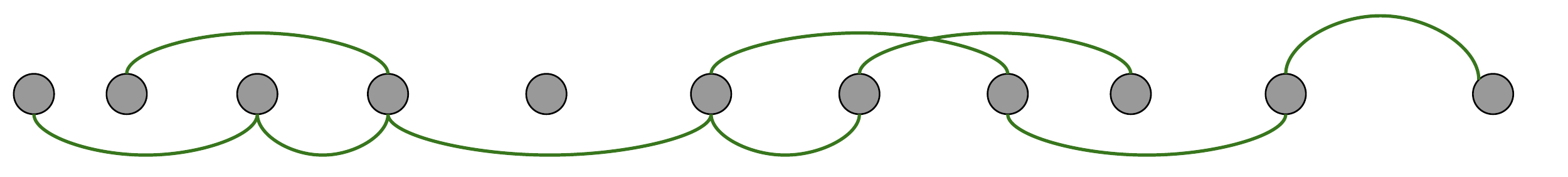}
    \caption{Isolated Node}
    \end{subfigure}
    
    \begin{subfigure}{0.45\textwidth}
    \includegraphics[width=\textwidth]{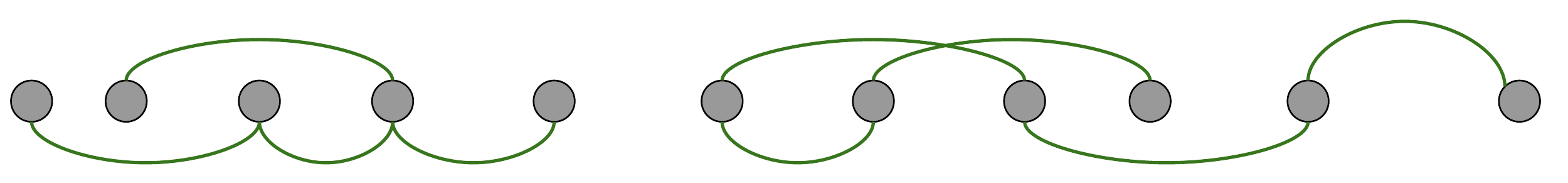}
    \caption{Uncrossed Gap}
    \end{subfigure}
    
    \begin{subfigure}{0.45\textwidth}
    \includegraphics[width=\textwidth]{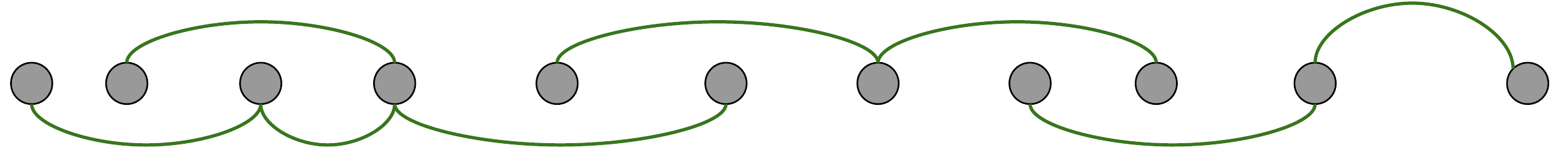}
    \caption{Split}
    \end{subfigure}
    \caption{Three different disconnection modes.}
    \label{fig:disconnection}
\end{figure}

Notice that in a hard RGG, if $v\in V$ is isolated and $u$ is the rightmost point of $V$ in $[0,v)$ (if there is one), then there is an uncrossed gap at $u$, and another one at $v$. Thus, a hard RGG may have uncrossed gaps without isolated nodes but not conversely; it is disconnected exactly when there is an uncrossed gap. The situation for soft RGGs is more complicated. It may be disconnected even if there are no isolated nodes or uncrossed gaps, as
illustrated in Figure \ref{fig:disconnection}c. Nevertheless, we conjecture that isolated nodes and uncrossed gaps together account for most of the probability of the graph being disconnected, converging to all of it in a suitable limiting regime. While we have been unable to prove this, we present some evidence below from simulations.

In order to gain an understanding of the causes of disconnection, we ran simulations for different system sizes and link ranges, for both Waxman and Rayleigh connection functions. We denote by $\piso$ the proportion of simulation runs (for a fixed set of parameter values) in which the graph contained an isolated node, by $\pucg$ the proportion with an uncrossed gap and by $\pisoucg$ the proportion with either. Comparing these with $\pfc$, the proportion of simulation runs in which the network fails to be fully connected, yields insights into the primary causes of disconnection.

\begin{figure*}
    \centering
    
    \begin{subfigure}{0.45\textwidth}
    \includegraphics[width=\textwidth]{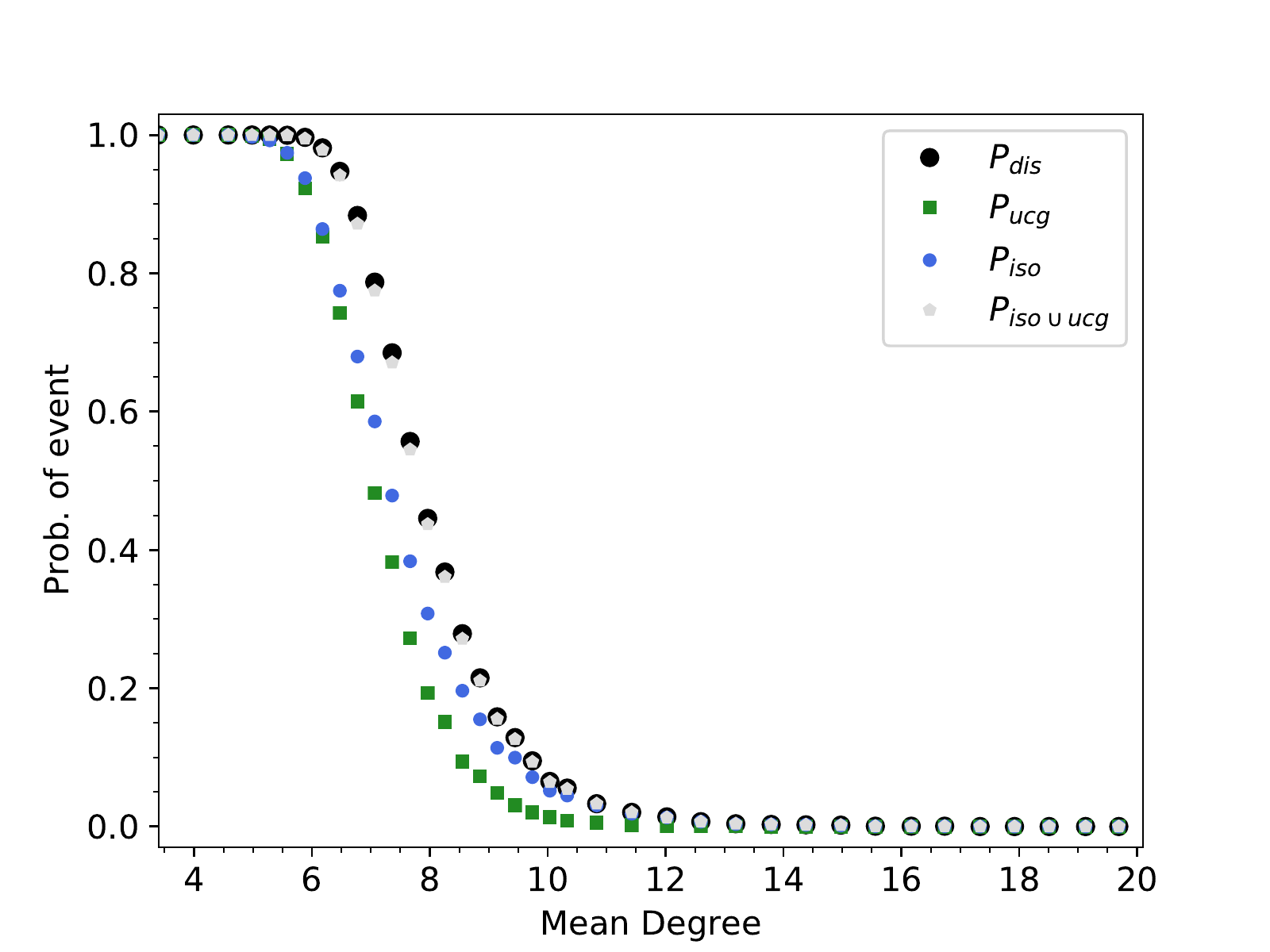}
    \caption{}
    \label{fig:WaxmanOverall_1000}
    \end{subfigure}
    \begin{subfigure}{0.45\textwidth}
    \includegraphics[width=\textwidth]{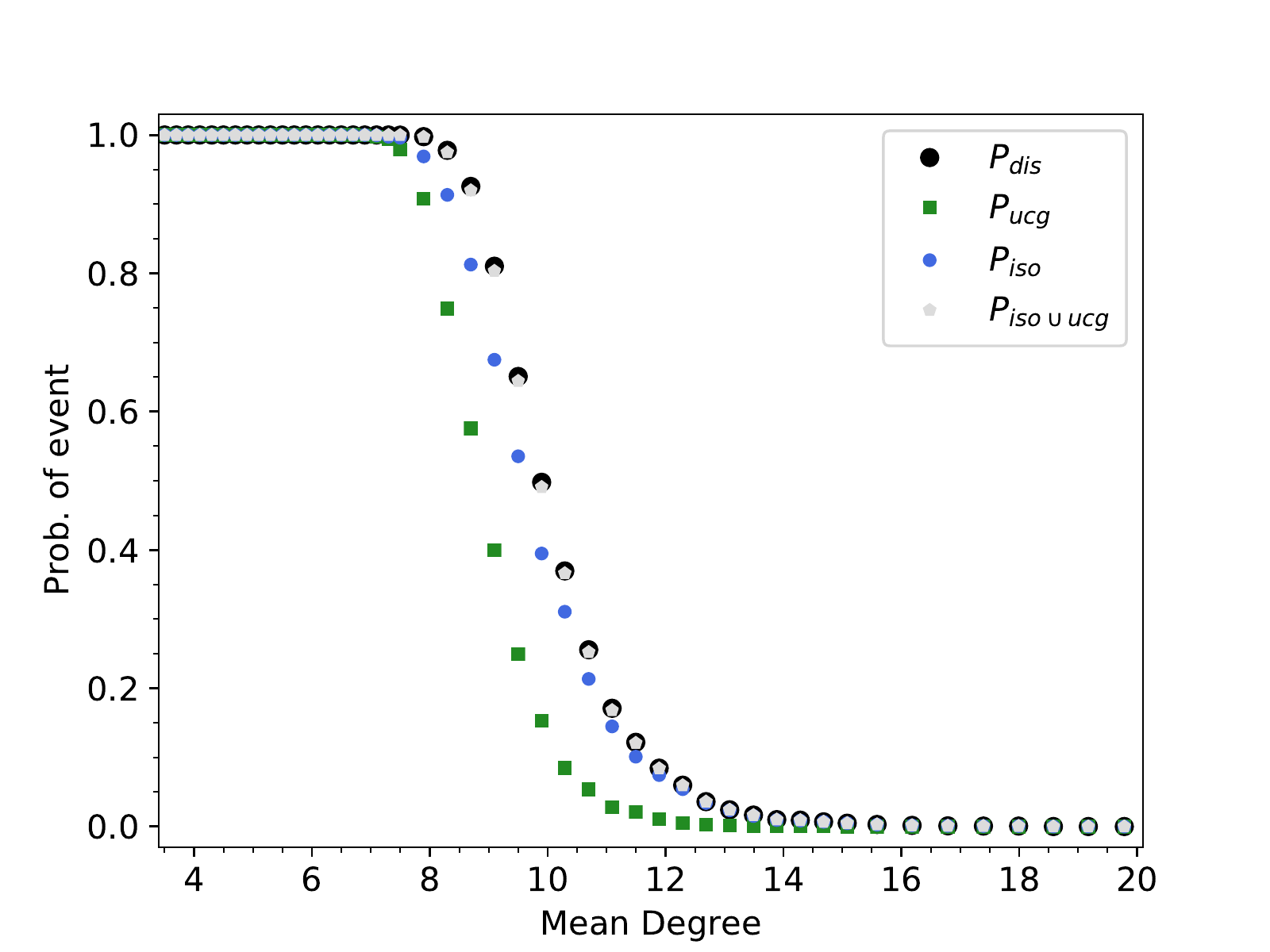}
    \caption{}
    \label{fig:WaxmanOverall_10000}
    \end{subfigure} \\
    \begin{subfigure}{0.45\textwidth}
    \includegraphics[width=\textwidth]{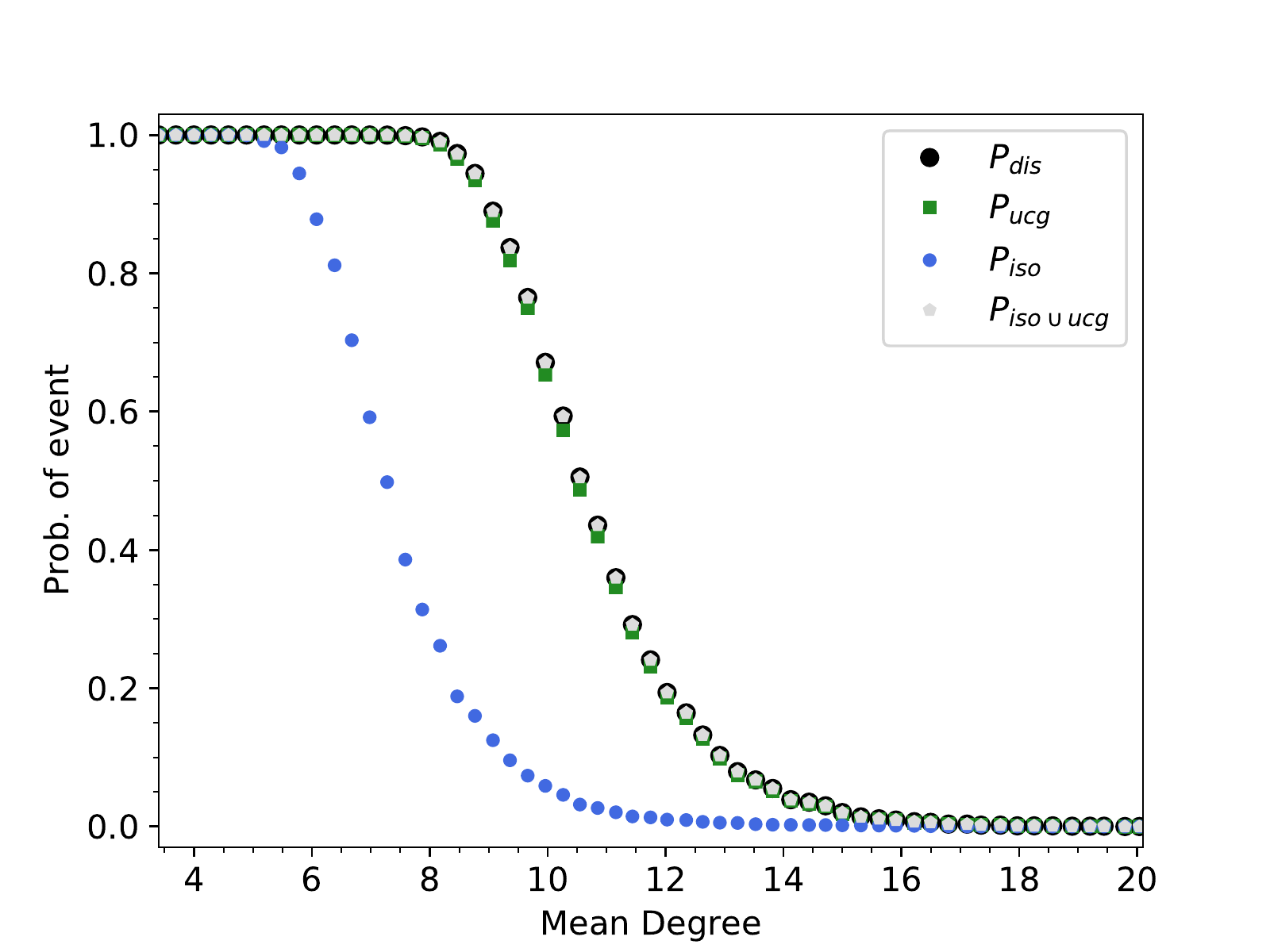}
    \caption{}
    \label{fig:RayleighOverall_1000}
    \end{subfigure}
    \begin{subfigure}{0.45\textwidth}
    \includegraphics[width=\textwidth]{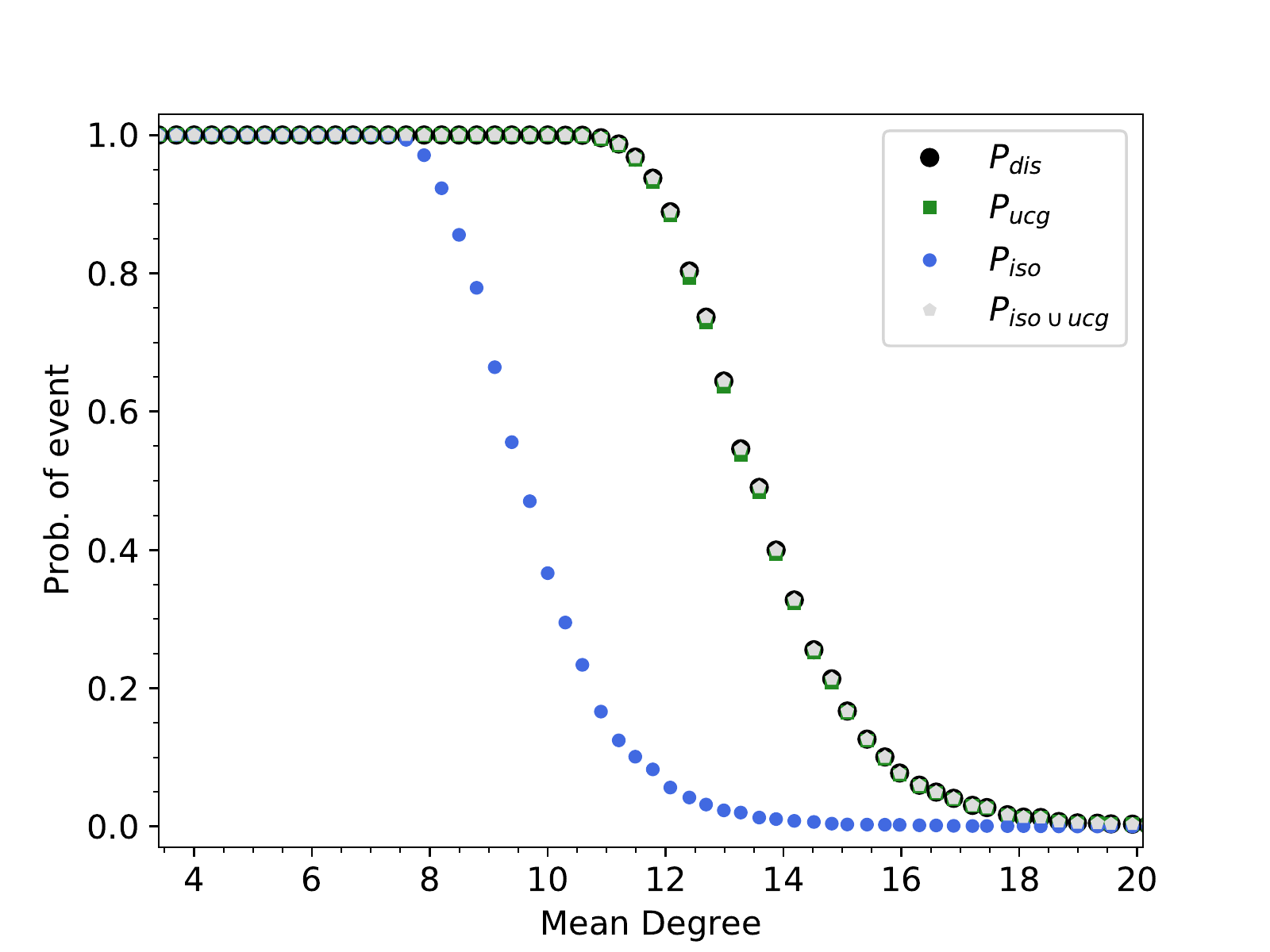}
    \caption{}
    \label{fig:RayleighOverall_10000}
    \end{subfigure}
    
    \caption{The proportion of simulations in which the network is disconnected ($\pfc$), compared with the proportion of simulations suffering from isolated nodes ($\piso$), uncrossed gaps ($\pucg$), or either of the two ($\pisoucg$). Each plot is based on $\nsims$ simulations of a soft RGG whose nodes are generated according to a unit rate PPP on a line segment of length $L = 1\,000$ or $L = 10\,000$, and whose edges are created according to a Waxman or Rayleigh connection function ($\eqn$~(\ref{conn_fun_wax_ray})). With (a) $L = 1000$, Waxman conn. fun., (b) $L = 10000$, Waxman conn. fun., (c) $L = 1000$, Rayleigh conn. fun., and (d) $L = 10000$, Rayleigh conn. fun.}
    \label{WaxmanOverall}
\end{figure*}

The simulations were run for two different system sizes, $L=1\,000$ and $L=10\,000$; since the Poisson process of node locations has unit intensity, $L$ is the expected number of nodes. The link range $r_c$ is varied, and the proportions of simulations exhibiting isolated nodes, uncrossed gaps or disconnection are plotted against the mean node degree (i.e. the average number of edges per node) corresponding to that value of $r_c$. (We chose not to plot the results against $r_c$ as the meaning of this parameter is somewhat opaque, whereas the mean node degree is intuitive.) The findings are displayed in Figure \ref{WaxmanOverall} for Waxman and Rayleigh connection functions as described in $\eqn$ (\ref{conn_fun_wax_ray}). The plots are based on $\nsims$ simulations for each set of parameter values. In both figures, we observe a fairly sharp transition between mean degrees for which disconnection (and its individual causes) have probability close to 1 and those for which it has probability close to 0. For isolated nodes, this transition occurs close to a mean degree of $\ln L$, exactly as for Erd\H{o}s-R\'enyi random graphs.



Figure \ref{fig:WaxmanOverall_1000} shows that in a system with $1\,000$ nodes on average and a Waxman connection function, isolated nodes and uncrossed gaps are approximately equally prevalent. But as the system size increases to $10\,000$ nodes on average, Figure \ref{fig:WaxmanOverall_10000} shows that isolated nodes become dominant. 
We also see that $\pisoucg$ is almost exactly the same as $\pfc$, supporting our intuition that node isolation and uncrossed gaps are the dominant modes of disconnection.

Figures \ref{fig:RayleighOverall_1000} and \ref{fig:RayleighOverall_10000} show that, for a Rayleigh connection function, the main obstruction to full connectivity appears to be uncrossed gaps. As the hard RGG is obtained in the limit as $\eta$ tends to infinity, and uncrossed gaps are the cause of disconnection in hard RGGs, it seems intuitive that they should dominate for large enough values of $\eta$. This intuition is supported by the simulations. Nevertheless, it is at odds with the analysis presented in Section \ref{sec:ucg}, which shows that isolated nodes dominate asymptotically whenever the connection function, $H(r)$, has unbounded support, i.e. it has infinite range. It appears then that the asymptotics don't kick in even at the large system sizes simulated. This includes simulations not shown in this paper looking at networks with up to $100\,000$ nodes. This discrepancy is discussed in further detail in Sections \ref{sec:ucg} and \ref{sec:conclusions} with some more detailed calculations given in the Appendices. Very large system sizes appear to be necessary to see this behaviour and hence it remains an open question to resolve this apparent discrepancy.

Once again, comparing $\pisoucg$ with $\pfc$ shows that the main contribution to disconnection comes from isolated nodes and uncrossed gaps. The remainder of the paper will focus on a theoretical analysis of these two events.

%% file: Section4_Iso.tex
\section{Isolated Nodes}
\label{sec:iso}

In this section, we derive bounds on the probability that there are isolated nodes, namely, nodes which have no edge to any other node.
We then explore an asymptotic regime in which $L$, the expected number of nodes in the system, tends to infinity, while the range of the connection function grows as a function of $L$. We seek to identify a scaling regime for the connection function at which a sharp threshold for the existence of isolated nodes can be discerned.



In order to eliminate boundary effects and simplify calculations, we modify our model slightly by identifying the endpoints of the line segment $[0,L]$, thereby enforcing periodic boundary conditions (PBCs). This is equivalent to modifying the line segment to be a 1-dimensional torus (or a circle). The effect of this modification on the number of isolated nodes is negligible, as it only affects nodes close to the boundary. We do not formalise this assertion here as this would add significant complexity but without much insight.

Denote by $\rho(x,y)$ the circular (or toroidal) distance between $x$ and $y$, and by $\rho_{\infty}(x,y)=|x-y|$ the Euclidean distance between them. Define 
$$
h=H\circ \rho,
$$
where $\circ$ denotes composition of functions; $h(x,y)$ is the probability that two points at locations $x$ and $y$ on the line segment with PBCs are connected. Note that $h(x,y)$ is symmetric, and invariant to translations on this line segment, i.e. $h(x,y) = h(y,x)$ and $h(x,y) = h(x-\delta, y-\delta)$ for some translation $\delta$. Denote by $\prob_x$ the Palm measure conditional on the Poisson process having a point at $x$, i.e., $\prob_x(E)$ is the probability of an event $E$ conditional on there being a node at $x$. The Palm measures conditional on any finite set of points are denoted analogously. We write $\E_{x}, \E_{x,y}$ etc. for the corresponding expectations. We make the following assumption about the connection function throughout this section.

\noindent{\bf Assumption A}: The function $H:\Reals_+ \mapsto [0,1]$ is integrable, i.e., 
$$
\|H\|_1 = \int_0^\infty H(x)dx < \infty.
$$

We now derive expressions for the expectation and variance of the number of isolated nodes. We will use these, along with Markov's and Chebyshev's inequalities, to obtain probability bounds on the existence of isolated nodes.


We start with the expected number of isolated nodes.
Let $\chi(x)$ denote the indicator that there is a point of the Poisson process at $x$ and that it is isolated. Now, conditional on there being a point at $x$, the remaining points constitute a unit rate Poisson process on $[0,L]$ by Slivnyak's theorem~\cite[Theorem 8.10]{haenggi2012stochastic}. Consequently, the set of neighbours of $x$ (the points to which an edge is present) constitute an inhomogenous Poisson process, with intensity measure $H(\rho(x,\cdot))=h(x,\cdot)$. The event that the point at $x$ is isolated is the event that this Poisson process is empty. Recall that for a Poisson process $X_{\nu}$ with intensity $\nu(\cdot)$ on $[0,L]$, the probability that this process is empty is given by 
$
\mathbb{P}(X_{\nu} = \emptyset) = e^{-\int_0^L \nu(x) dx}.
$
The probability that the node at $x$ is isolated is the expectation of the indicator of this event:
\begin{equation} \label{eq:isol_prob}
\begin{aligned}
\E_{x}[\chi(x)] &=\mathbb{P}_x(\chi(x) = 1) \\ &=e^{-\int_0^L h(x,y)dy} = e^{-\int_0^L h(0,y)dy},
\end{aligned}
\end{equation}
which does not depend on $x$, as expected by the symmetry in the model. 

Let $\niso$ denote the number of isolated nodes. By integrating $\eqn$ (\ref{eq:isol_prob}) over $x\in [0,L]$ with respect to the intensity of the Poisson process, which was assumed to be unity, we get
\begin{equation}
    \label{eqn:E_N_iso}
    \mathbb{E}[\niso] = \int_0^L \mathbb{E}_x(\chi(x)) dx = L e^{-2\int_0^{L/2} H(y)dy},
\end{equation}
where the last equality holds because we are working with toroidal distance. As $L$ tends to infinity, the exponent tends to $-2\| H\|_1$, which is a finite constant. Hence, $\E[\niso]$ scales linearly in $L$. This leads us to ask if the random variable $\niso$ does so as well, i.e., if it behaves as an extensive quantity, in the language of thermodynamics. The answer depends on whether correlations decay quickly enough, and can be addressed by calculating the variance of $\niso$.

As a first step towards computing the variance, we condition on there being points at $x$ and $y$; by Slivnyak's theorem again, the Palm measure $\prob_{x,y}$ corresponds to a unit rate Poisson process on $[0,L]\setminus \{ x,y\}$. The set of these points that are neighbours of either $x$ or $y$ constitute an inhomogenous Poisson process of intensity $h(x,\cdot)+h(y,\cdot)-h(x,\cdot)h(y,\cdot)$. The event that the points at both $x$ and $y$ are isolated is the event that this point process is empty, and that there is no edge between $x$ and $y$. Hence, 
\begin{equation}\label{eq:isol_pair}
\begin{aligned}
&\E_{x,y}[\chi(x)\chi(y)] = \\
&= (1-h(x,y)) e^{-\int_0^L [h(x,z)+h(y,z)-h(x,z)h(y,z)]dz} \\
&\leq e^{-\int_0^L [h(x,z)+h(y,z)-h(x,z)h(y,z)]dz},
\end{aligned}
\end{equation}
since $0\leq h(x,y) \leq 1$. Then, using translation invariance,
\begin{equation} \label{eq:niso_stats}
\begin{aligned}
\E[&\niso(\niso-1)] \\
&= \int_0^L \int_0^L \E_{x,y}(\chi(x),\chi(y)) dxdy \\
&= L \int_0^L \E_{0,x}[\chi(0)\chi(x)]dx \\
&\leq L \int_0^L e^{-\int_0^L [h(0,z)+h(x,z)-h(0,z)h(x,z)]dz} dx.
\end{aligned}
\end{equation}

Recall that the squared coefficient of variation of a random variable, denoted $c_V^2$, is defined as the ratio of its variance to its squared mean. Hence, we obtain from $\eqns$ (\ref{eq:isol_prob}), (\ref{eqn:E_N_iso}), (\ref{eq:isol_pair}), (\ref{eq:niso_stats}) that 
\begin{equation} \label{eq:niso_cvsq}
\begin{aligned}
&c_V^2(\niso) - \frac{1}{\E[\niso]}= \\
&= \frac{1}{L}\int_0^L \frac{\E_{0,x}[\chi(0)\chi(x)]}{\E_0[\chi(0)]^2}dx - 1\\
&\leq \frac{1}{L}\int_0^L \left( \frac{e^{-\int_0^L [h(0,z))+h(x,z)-h(0,z)h(x,z)] dz}}{e^{-2\int_0^L h(0,z)dz}} - 1 \right)dx \\
&= \frac{1}{L}\int_0^L \left( e^{\int_0^L h(0,z)h(x,z)dz} -1 \right)dx.
\end{aligned}
\end{equation}
Since $h$ is bounded above by 1, we have
\begin{equation*}
\begin{aligned}
\int_0^L h(0,z)h(x,z)dz & \leq \int_0^L h(0,z)dz  \\
& \leq 2\int_0^{\infty} H(z)dz < \infty.
\end{aligned}
\end{equation*}
In other words, the above integral is bounded, uniformly in $x$. Hence, there is a finite constant $C$, which does not depend on $x$ or $L$, such that
$$
e^{\int_0^L h(0,z)h(x,z)dz} \leq 1+C\int_0^L h(0,z)h(x,z)dz.
$$
It follows that 
\begin{align*}
    \frac{1}{L}\int_0^L & \Bigl( e^{\int_0^L h(0,z)h(x,z)dz} -1 \Bigr)dx \\
    &\leq \frac{1}{L}\int_0^L \left( C\int_0^L h(0,z)h(x,z)dz \right)dx \\
    &= \frac{C}{L}\int_0^L h(0,z) \left( \int_0^L h(x,z) dx \right) dz \\
    &= \frac{4C}{L} \int_0^{L/2} H(z)dz \int_0^{L/2} H(y)dy \\
    & \leq \frac{4C}{L} \| H\|_1^2.
\end{align*}
The interchange of the order of integration in the third line is justified by Tonelli's theorem. Substituting the above in $\eqn$ (\ref{eq:niso_cvsq}), we obtain that 
\begin{equation} \label{eq:niso_cvsq_bd}
c_V^2(\niso) \leq \frac{1}{\E[\niso]} + \frac{4C}{L} \| H\|_1^2.
\end{equation}
We already noted that $\E[\niso]$ scales linearly in $L$. Therefore, the right hand side above tends to zero as $L$ tends to infinity, Now, by Chebyshev's inequality, we have for any $\epsilon>0$ that 
\begin{equation*}
\prob\left(|\niso-\E [\niso]|> \epsilon \E [\niso]\right) \leq \frac{c_V^2(\niso)}{\epsilon^2}, 
\end{equation*}
which tends to zero as $L$ tends to infinity. Thus, $\niso$ concentrates around its expected value, and scales linearly with $L$, making it an extensive quantity.

We now turn to the question of identifying a scaling regime where the probability of seeing isolated nodes in the graph exhibits a sharp transition.

\subsection{Scaled Connection Function}

Consider a family of 1-D soft RGGs as above, indexed by $L$, and with scaled connection functions $H^L(z)=H(z/R_L)$, where $R_L$ is an increasing function of $L$. For now, the only assumption we make is that $R_L$ tends to infinity and $R_L/L$ tends to zero as $L$ tends to infinity. Let $h^L=H^L \circ \rho$. We shall study these graphs in the asymptotic regime $L\to \infty$.

We begin by rewriting $\eqn$ (\ref{eqn:E_N_iso}) as
\begin{equation*}
\begin{aligned}
    \E[\niso]&= L\exp \left( -2\int_0^{L/2} H^L(y)dy \right)\\
    & =L\exp \left( -2\int_0^{L/2} H(y/R_L)dy \right)
\end{aligned}
\end{equation*}
Making the change of variables $u = y/R_L$, we get
\begin{equation}
\label{eqn:E_niso_scaled}
\mathbb{E}[\niso] = L\exp\left(-2 R_L \int_0^{L / 2R_L} H(u) du \right).
\end{equation}
Since $R_L/L$ was assumed to tend to zero, $\int_0^{L / 2R_L} H(u) du$ tends to $\|H\|_1$ as $L$ tends to infinity, and we get
\begin{equation*}
\lim_{L \rightarrow \infty} \mathbb{E}[\niso] = \lim_{L \rightarrow \infty}Le^{-2 R_L \|H\|_1},
\end{equation*}
provided the limit on the right exists. 

We see from the above expression that the scaling required is $R_L$ growing logarithmically in $L$. Indeed, taking $R_L=\gamma \ln L$, we get
\begin{equation}
\label{eqn:critical_expectation}
\E[\niso]  \xrightarrow{L\rightarrow \infty}
	\begin{cases}
		0, &\mbox{ if } 2\gamma \|H\|_1 > 1, \\
		1, &\mbox{ if } 2\gamma \|H\|_1 = 1, \\ 
		\infty, &\mbox{ if } 2\gamma \|H\|_1 <1.
	\end{cases}
\end{equation}
Thus, the mean number of isolated nodes exhibits a sharp transition at $\gamma=1/2\| H\|_1$. We wish to show that the random variable denoting the number of isolated nodes does so as well.

It is immediate from Markov's inequality, which states that $\mathbb{P}(X \geq 1) \leq \mathbb{E}[X]$ for any non-negative random variable $X$, that
\begin{equation} \label{eq:markov_bd}
\mathbb{P}(\niso\geq 1) \xrightarrow{L \rightarrow \infty} 0 \enspace \mbox{ if } \enspace 2\gamma \| H\|_1>1.
\end{equation}
In words, such a choice of $\gamma$ ensures that large networks have a vanishingly small chance of containing isolated nodes.

To investigate the behaviour for $\gamma$ smaller than this threshold, we must return to the coefficient of variation. We can rewrite $\eqn$ (\ref{eq:niso_cvsq}) as 
\begin{equation} \label{eq:niso_cvsq_crit}
\begin{aligned}
c^2_V&(\niso) - \frac{1}{\E[\niso]} \\
&\leq \frac{1}{L} \int_0^L \Bigl( e^{\int_0^L h^L(0,z)h^L(x,z)dz}-1 \Bigr)dx.
\end{aligned}
\end{equation}
Defining
\begin{equation*} 
g^L(x) = \int_0^L h^L(0,z)h^L(x,z)dz,
\end{equation*}
we have
\begin{equation} \label{eq:scaled_convol_bd}
\begin{aligned}
    \int_0^L g^L(x)dx &= \int_0^L \int_0^L h^L(0,z)h^L(x,z)dzdx \\
    &= \int_0^L h^L(0,z) \left( \int_0^L h^L(x,z)dx \right)dz \\
    &= \left( \int_0^L h^L(0,y)dy \right)^2 \\
    &= \left(2\int_0^{L/2} H\left(\frac{y}{R_L}\right)dy\right)^2 \\
    &= \left(2R_L \int_0^{L/2R_L} H(y)dy \right)^2 \\
    & \leq 4R_L^2 \|H\|_1^2.
\end{aligned}
\end{equation}
The interchange of the order of integration in the second line is justified by Tonelli's theorem.
For $\alpha>0$, define the set $\Psi_{\alpha}= \{ x\in [0,L):g^L(x)> \alpha \}$. Since $g^L$ is non-negative, it follows from $\eqn$ (\ref{eq:scaled_convol_bd}) that the Lebesgue measure of the set $\Psi_{\alpha}$, denoted $m(\Psi_{\alpha})$, is bounded as follows:
\begin{equation} \label{eq:tail_bound_1}
m(\Psi_{\alpha}) \leq \frac{4 R_L^2 \| H\|_1^2}{\alpha}.
\end{equation}
Now, if $x$ is not in $\Psi_{\alpha}$, then $g^L(x) \leq \alpha$, and so there is a constant $C_{\alpha}>0$ such that $e^{g^L(x)} \leq 1+C_{\alpha}g^L(x)$. Hence,
\begin{equation} 
\label{eq:int_psicomp_bd}
\begin{aligned}
\int_{\Psi_{\alpha}^c} \Bigl( e^{g^L(x)}-1 \Bigr)dx &\leq C_{\alpha}\int_{\Psi_{\alpha}^c} g^L(x)dx \\
&\leq 4C_{\alpha}R_L^2 \|H\|_1^2,
\end{aligned}
\end{equation}
where $\Psi_{\alpha}^c$ denotes the complement of $\Psi_{\alpha}$ in $[0,L]$. We have used $\eqn$ (\ref{eq:scaled_convol_bd}) and the non-negativity of $g^L$ to obtain the last inequality.

Integrating $e^{g^L(x)} - 1$ over $\Psi_\alpha$, we have
\begin{equation}
    \label{eq:bound_over_psi_alpha}
    \begin{aligned}
        \int_{\Psi_\alpha} & \Bigl( e^{g^L(x)} - 1 \Bigr) dx  \\
        & = \int_{\Psi_\alpha} \Bigl( \exp\Bigl(\int_0^L h^L(0,z) h^L(x,z) \Bigr) - 1 \Bigr) dx \\
        & \leq \int_{\Psi_\alpha} \Bigl( \exp\Bigl(\int_0^L h^L(0,z) \Bigr) - 1 \Bigr) dx \\
        & = \int_{\Psi_\alpha} \Bigl( \exp\Bigl(2R_L\int_0^{L/2R_L} H(z) dz \Bigr) - 1 \Bigr) dx \\
        & \leq m(\Psi_\alpha) \bigl( e^{2R_L \|H_1\|} - 1 \bigr) \\
        & \leq \frac{4 R_L^2 \| H\|_1^2}{\alpha} \bigl( e^{2R_L \|H_1\|} - 1 \bigr).
    \end{aligned}
\end{equation}
We have used the fact that $h^L(\cdot,\cdot) \leq 1$ to obtain the inequality on the third line, and $\eqn$ (\ref{eq:tail_bound_1}) to obtain the last inequality. Now, substituting $\eqns$ (\ref{eq:int_psicomp_bd}) and (\ref{eq:bound_over_psi_alpha}) into eqn.  (\ref{eq:niso_cvsq_crit}), and noting that $\Psi_{\alpha} \cup \Psi_{\alpha}^c=[0,L]$, we get
\begin{equation*}
    \begin{aligned}
    c^2_V(\niso) \leq  & \frac{1}{\E[\niso]} + \frac{4C_{\alpha} R_L^2 \|H\|_1^2}{L}\\
    & \quad +\frac{4 R_L^2 \| H\|_1^2 }{\alpha L} \bigl( e^{2R_L \|H_1\|} - 1 \bigr) .
    \end{aligned}
\end{equation*}
Thus, for the scaling $R_L=\gamma \ln L$, we have
\begin{equation*}
\begin{aligned}
c^2_V(\niso) \leq &  \frac{1}{\E[\niso]} + \frac{4C_{\alpha} \gamma^2 \|H\|_1^2 \ln^2 L}{L} \\
& \qquad + \frac{4 \gamma^2 \| H\|_1^2 \ln^2 L}{\alpha L} \bigl( e^{2\gamma \|H_1\| \ln L} - 1 \bigr).
\end{aligned}
\end{equation*}

Suppose $2\gamma \| H\|_1<1$. Then, as $L$ tends to infinity, $\E[\niso]$ tends to infinity by $\eqn$  (\ref{eqn:critical_expectation}),
while the second and third terms in the sum on the RHS tend to zero. Hence, 
\begin{equation} \label{eq:cv_final_bound}
c^2_V(\niso) \xrightarrow{L\rightarrow \infty} 0, \mbox{ if } 2\gamma\| H\|_1<1.
\end{equation}
But, by Chebyshev's inequality,
\begin{equation*}
\prob(\niso =0 ) \leq \prob(|\niso-\E \niso|\geq \E \niso) \leq c_V^2(\niso),
\end{equation*} 
and so, it tends to zero as $L$ tends to infinity, for $\gamma<1/2\| H\|_1$. Combining this with the result established in $\eqn$ (\ref{eq:markov_bd}), we conclude that
\begin{equation}
\label{eqn:P_iso_scaled_limit}
\prob(\niso=0) \xrightarrow{L \rightarrow \infty}
\begin{cases}
0, & \mbox{if } 2\gamma\| H\|_1<1, \\
1, & \mbox{if } 2\gamma\| H\|_1>1.
\end{cases}
\end{equation}

\begin{figure*}
    \centering
    \begin{subfigure}{0.45\textwidth}
    \includegraphics[width=\textwidth]{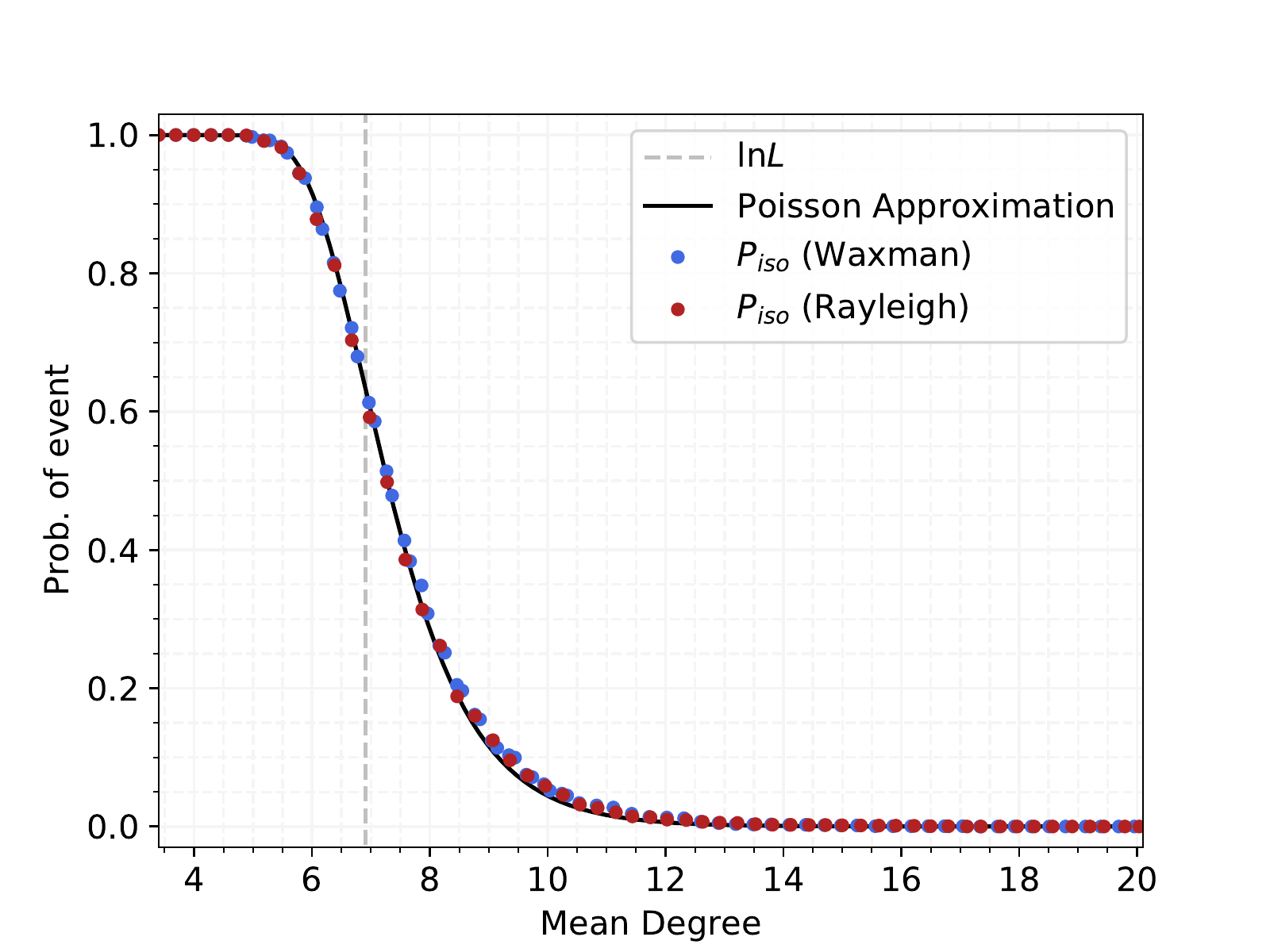}
    \caption{}
    \end{subfigure}
    \begin{subfigure}{0.45\textwidth}
    \includegraphics[width=\textwidth]{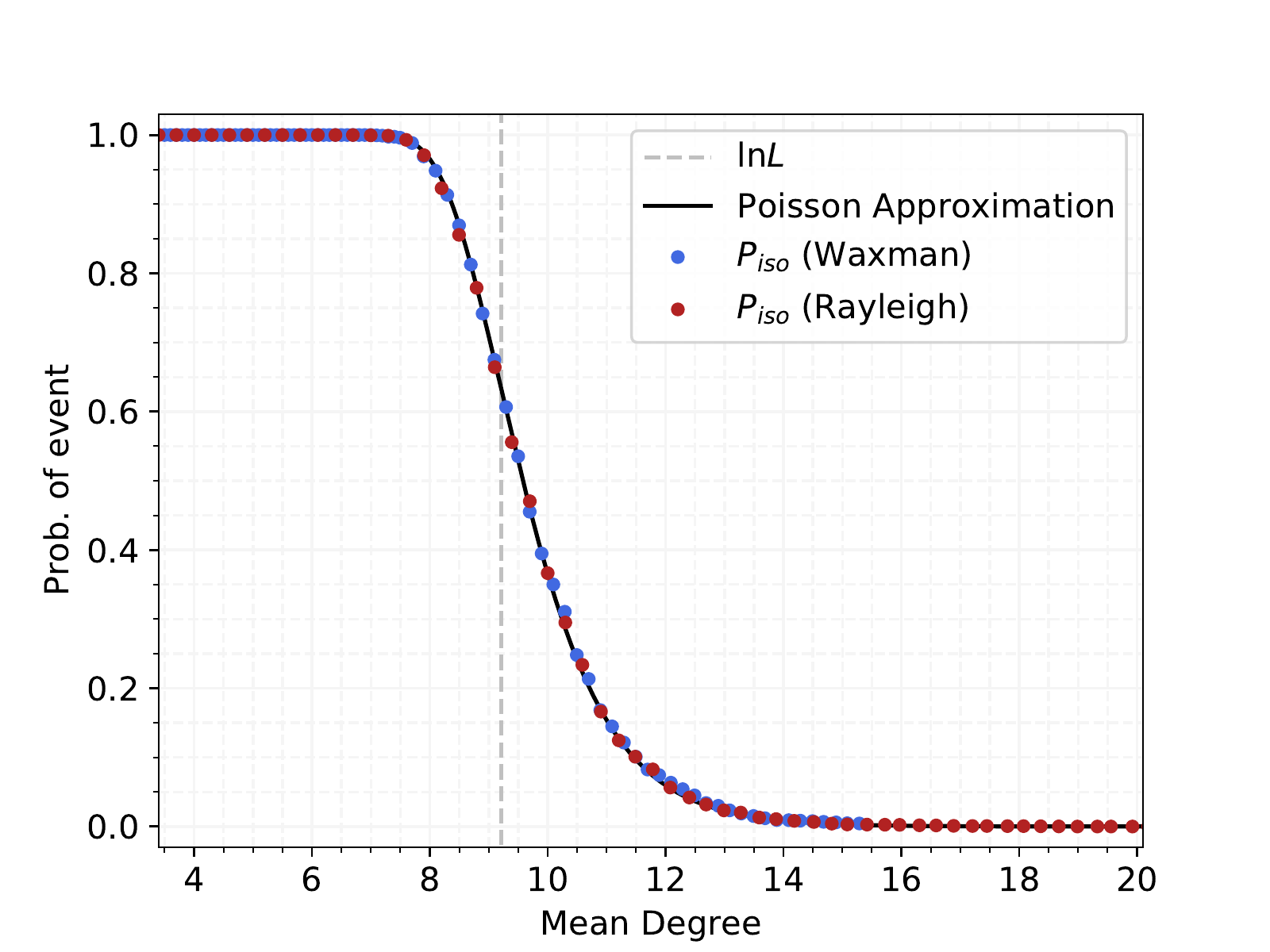}
    \caption{}
    \end{subfigure}
    \caption{Probability of isolated nodes being present for line segments of length (a) $1\,000$ and (b) $10\,000$: comparison of simulations (circles) with Poisson approximation (black solid line). Each circle is a proportion from $\nsims$ simulation runs.}
    \label{fig:iso_wax}
\end{figure*}

\medskip
\noindent{\bf Remarks.} 
\begin{enumerate}
    \item
    The mean degree of the nodes, given by $\int_0^L h^L(0,z)dz$, is asymptotic to $2\gamma \|H\|_1 \ln L$, while the mean number of nodes in the interval $[0,L]$ is equal to $L$. Thus, the theorem states that there is a sharp threshold for the existence of isolated nodes when the mean degree is equal to the natural logarithm of the mean number of nodes. This is shown in Figures \ref{fig:iso_wax}(a) and (b) by the dashed vertical line.  This threshold is exactly the same as for the hard RGG (see, e.g, \cite{appel2002connectivity, han2006very}).
    \item
    As was stated above, a superficially different scaling is seen in the densification regime. Here, the line segment remains fixed as $[0,1]$, whilst the intensity $\lambda$ of the PPP increases to infinity, and the scaling parameter $R(\lambda)$ decreases to zero. The results given here can be translated to this regime using a suitable rescaling of $H(.)$. See \cite{mao2011towards} for a more detailed discussion.
    \item
    The threshold for the emergence of isolated nodes is insensitive to the connection function $H$ and depends only on its integral, which is the expected number of connections. This is demonstrated in Figures \ref{fig:iso_wax}(a) and (b) by the fact that the simulations for both the Waxman and Rayleigh connection functions overlap with each other. This lack of sensitivity is also seen in the threshold for soft RGGs in two dimensions which depends only on the integral of $H$ as was shown in \cite{mao2011towards}.
\end{enumerate}

The threshold phenomenon arises in the limit of large system sizes (with the connection function scaling suitably). It does not provide explicit probabilities for the occurrence of isolated nodes in finite-sized networks. Moreover, the bounds provided by Markov's and Chebyshev's inequalities are rather weak. Finally, the limit result does not say what happens at $\gamma=1/2\| H\|_1$.

In \cite{penrose2016connectivity}, it was shown for soft RGGs in two or more dimensions, and for a large class of connection functions, that the number of isolated nodes can be well approximated by a Poisson random variable; in particular, the probability that there are no isolated nodes is well approximated by $e^{-\lambda}$, where $\lambda$ is the mean of this Poisson distribution. We conjecture that the same is true in our model, and provide a precise statement below.


\begin{conjecture}
\label{prop:IsoNodesPoisson}
Fix $\tau \in \Reals_+$ and consider a soft RGG on $[0,L]$, with connection functions $H(\cdot/R_L)$, where $R_L=\ln(\tau L)/2\| H\|_1$. Let $\niso^L$ denote the number of isolated nodes, whose dependence on $L$ has been made explicit in the notation. Then, for any sequence $L$ tending to infinity, the sequence $\niso^L$ converges in distribution to a Poisson distribution with mean $1/\tau$. In particular, $\prob(\niso^L=0)$ tends to $e^{-1/\tau}$.
\end{conjecture}

Evidence in support of this conjecture is presented in Figure \ref{fig:iso_wax}, where we have compared the proportion of simulations with no isolated nodes with this Poisson approximation for systems with $L=1\,000$ and $L=10\,000$. The proportion of $\nsims$ simulations in which at least one isolated node was present is plotted against the mean node degree, for both the Waxman (blue circles) and Rayleigh (red circles) connection functions. The plots demonstrate that the probability of isolated nodes appearing in the network is insensitive to whether the connection function is Waxman or Rayleigh, and depends only on the mean degree. The plots also show that the probability of occurrence of isolated nodes rises sharply as the mean degree falls below $\ln L$, and that the transition becomes sharper as the system size grows. Finally, the solid black line in the plots depicts the function $1 - \exp\left( -L e^{-\bar{k}} \right)$, where $\bar{k}$ is the mean degree. This function is the Poisson approximation for the probability of at least one isolated node being present, as stated in Conjecture~\ref{prop:IsoNodesPoisson}. The plots show that the Poisson approximation is very close to the observed prevalence of isolated nodes in the simulations.

%% file: Section6_UCG_New.tex
\section{Uncrossed Gaps}
\label{sec:ucg}

In the previous section, we showed that the probability that isolated nodes exist exhibits a sharp transition, from being close to 1 to being close to 0, as the mean degree is increased. Moreover, the point at which this transition occurs depends only on the mean degree, and is insensitive to the connection function.

Our primary interest is in the probability of disconnection, not just in whether isolated nodes exist. We noted in Section~\ref{sec:analysis} that there are other modes of disconnection (Figure~\ref{fig:disconnection}), and conjectured that uncrossed gaps are the other main mode. Hence, we now turn to estimating the probability that there are uncrossed gaps. We will focus on the scaling regime identified in the last section, namely $R_L=\gamma \ln L$, and $H^L(\cdot)=H(\cdot/R_L)$.
In this section, we need to assume the following about the connection function.

\noindent{\bf Assumption B} The connection function $H:[0,\infty)\to [0,1]$ is monotone decreasing and has unbounded support, i.e., $\int_r^{\infty} H(x)dx>0$ for all $r\geq 0$.

Our main result in this section is that the probability of seeing uncrossed gaps becomes negligible as $L$ tends to infinity, for $R_L=\gamma \ln L$ and any $\gamma>0$. In particular, it holds for $\gamma= 2\| H\|_1$, when isolated nodes begin to emerge. Consequently, the main cause of disconnection of soft RGGs in 1-D is isolated nodes, just as in higher dimensions; 1-D hard RGGs are anomalous in this regard.

We remark that, while we believe the assumption about unboundedness of the support of $H$ to be essential (indeed, our result does not hold for hard RGGs), the assumption of monotonicity appears to be an artefact of our proof technique.

In order to avoid boundary effects near 0 and $L$, we extend the Poisson process from $[0,L]$ to the infinite real line. (In this case, it is less convenient to work with the circle, since the existence of a single uncrossed gap on the circle does not guarantee disconnection.) We want to calculate the probability that there are no edges crossing the origin, under the Palm measure corresponding to having a point of the Poisson process at the origin. We denote probabilities and expectations under this measure by $\prob_0$ and $\E_0$ respectively, as before.

Let $\Phi$ denote a Poisson point process on $\Reals$ of unit intensity. For $A\subset \Reals$, we write $\Phi_A$ to denote the restrictions of $\Phi$ to $A$. Denote by $X_0$ the indicator of the event that there is a point at the origin and that it marks an uncrossed gap, i.e., there are no edges between $\Phi_{(-\infty,0)}\cup \{ 0\}$ and $\Phi_{(0,\infty)}$. We wish to compute $\E_0[X_0]$. As it is not amenable to exact calculation, we obtain a bound on it below.

Fix $\alpha, \delta>0$ and set $n_L=\lfloor \delta R_L \rfloor$. Denote by $X^{\alpha}_0$ the indicator of the event that there are no edges between $\Phi_{(-\alpha R_L,0)} \cup \{ 0\}$ and $\Phi_{(0,\infty)}$.  Clearly, $X_0\leq X^{\alpha}_0$. Moreover, conditional on $\Phi_{(-\alpha R_L,0)}$, the set of points on $(0,\infty)$ which have an edge to some point in $(-\alpha R_L, 0]$ constitute a Poisson process with intensity 
\begin{equation} \label{eq:connected_pts_intensity}
\begin{aligned}
\lambda_{\Phi_{(-\alpha R_L,0)}} (y) & = 1-\bar{h}^L(0,y) \prod_{z\in \Phi_{(-\alpha R_L,0)}}\bar{h}^L(z,y) \\
& = 1-\bar{H}^L(y) \prod_{z\in \Phi_{(-\alpha R_L,0)}}\bar{H}^L(y-z),
\end{aligned}
\end{equation}
where $\bar{h}^L(0,y) = 1 - h^L(0,y)$ and $\bar{H}^L(y) = 1 - H^L(y)$. Hence, the number of such points is a Poisson random variable with mean $\int_0^{\infty} \lambda_{\Phi_{(-\alpha R_L,0)}}(y)dy$. Since $X^{\alpha}_0$ is the indicator that this random variable takes the value zero, we obtain that 
\begin{equation} \label{eq:split_cond_bd1}
\E_0 \bigl[ X^{\alpha}_0 \bigm| \Phi_{(-\alpha R_L,0)} \bigr] = \exp \left( -\int_0^{\infty} \lambda_{\Phi_{(-\alpha R_L, 0)}}(y) dy \right).
\end{equation}
We also have by $\eqn$ (\ref{eq:connected_pts_intensity}) and Assumption B that 
\begin{equation} \label{eq:intensity_bd}
\lambda_{\Phi_{(-\alpha R_L,0)}}(y) \geq 1-\bar{H}^L(y+\alpha R_L)^{1+|\Phi_{(-\alpha R_L,0)}|},
\end{equation}
where $|\Phi_A|$ denotes the number of points of $\Phi$ in the set $A$. Hence, if we condition further on there being at least $n_L$ points of $\Phi$ in $(-\alpha R_L,0)$, then we obtain by substituting $\eqn$ (\ref{eq:intensity_bd}) in $\eqn$ (\ref{eq:split_cond_bd1}) that
\begin{equation} \label{eq:split_cond_bd2}
\begin{aligned}
\E_0 \bigl[ X^{\alpha}_0 \bigm| |\Phi_{(-\alpha R_L,0)}| \geq n_L \Bigr] &\leq \exp \Bigl( -\int_0^{\infty} \bigl( 1-\bar{H}^L(y+\alpha R_L)^{1+n_L} \bigr)dy \Bigr) \\
&= \exp \Bigl( -R_L\int_0^{\infty}  \bigl(1-\bar{H}(z+\alpha)^{1+n_L} \bigr)dz \Bigr).
\end{aligned}
\end{equation}
We have used the change of variables $z=y/R_L$ to obtain the last equality.

Define $H^{-1}(x)=\sup \{z: H(z)\geq x \}$ to be the generalised inverse of $H$; $H^{-1}$ maps $(0,H(0)]$ to $[0,\infty)$. Since $R_L=\gamma \ln L$ and $n_L = \lfloor \delta R_L \rfloor$ for fixed constants $\gamma, \delta>0$, $n_L$ tends to infinity as $L$ tends to infinity. Consequently, 
\begin{equation}
\label{eq:Hinv_tail}
    H^{-1} \left(\frac{1}{1+n_L}\right) \to \infty \mbox{ as } L \to \infty, 
\end{equation}
because the support of $H$ is unbounded by Assumption B. Since $H$ is monotone decreasing, we have for all $y\leq H^{-1}\left(\frac{1}{1+n_L}\right)$ that $H(y)\geq \frac{1}{1+n_L}$, and hence, $\bar{H}(y)^{1+n_L}\leq e^{-1}$. Therefore, 
\begin{equation*}
\begin{aligned}
R_L \int_0^{\infty}  \bigl(1-\bar{H}(z+\alpha)^{1+n_L} \bigr)dz & \geq R_L \int_0^{H^{-1}(\frac{1}{1+n_L})}  \bigl(1-\bar{H}(z+\alpha)^{1+n_L} \bigr)dz \\
& \geq (1-e^{-1})\gamma H^{-1}\left( \frac{1}{1+n_L} \right) \ln L
 = \omega(\ln L),
\end{aligned}
\end{equation*}
where the last equality follows from $\eqn$ (\ref{eq:Hinv_tail}). We use the notation $f(x)=\omega(g(x))$ to denote, for functions $f, g: \Reals_+\to \Reals_+$, that $g(x)/f(x)$ tends to infinity as $x$ tends to infinity. Substituting the above expression into $\eqn$ (\ref{eq:split_cond_bd2}), we obtain that 
\begin{equation} \label{eq:split_cond_bd3}
\E_0 \bigl[ X^{\alpha}_0 \bigm| |\Phi_{(-\alpha R_L,0)}| \geq n_L \Bigr] = o(L^{-\kappa}),
\end{equation}
for arbitrary $\kappa>0$.

Next, we bound the probability of the event that $|\Phi_{(-\alpha R_L,0)}|$ is smaller than $n_L$. Since $|\Phi_{(-\alpha R_L,)}|$ is a Poisson random variable with mean $\alpha R_L$, a standard use of the Bernstein inequality (also known as Chernoff's inequality) yields that 
\begin{equation} \label{eq:poisson_tail_bd}
\begin{aligned}
\prob ( |\Phi_{(-\alpha R_L,0)}|< n_L ) &\leq \prob( |\Phi_{(-\alpha R_L,0)}| \leq \delta R_L) \\
&\leq \inf_{\theta\leq 0} e^{-\theta \delta R_L} \E \bigl[ e^{\theta |\Phi_{(-\alpha R_L,0)}|} \bigr] \\
&= \exp \bigl(- \sup_{\theta \leq 0} \, \bigl[ \theta \delta R_L -\alpha R_L(e^{\theta}-1) \bigr] \bigr) \\
&= \exp ( -\gamma I_{\alpha}(\delta) \ln L),
\end{aligned}
\end{equation}
where $I_{\alpha}(\delta) = -\delta \ln \frac{\alpha}{\delta} + \alpha - \delta$. It is easy to see that, for fixed $\alpha>0$, $I_{\alpha}(\delta)$ tends to $\alpha$ as $\delta$ decreases to zero. Hence, given arbitrary $\kappa, \gamma>0$, we can choose $\alpha$ and $\delta$ such that $\gamma I_{\alpha}(\delta)>\kappa$. Then,
\begin{equation} \label{eq:poisson_tail_bd2}
\prob( |\Phi_{(-\alpha R_L,0)}|<n_L ) \leq e^{-\kappa \ln L} = L^{-\kappa}. 
\end{equation}
Combining $\eqns$ (\ref{eq:split_cond_bd3}) and (\ref{eq:poisson_tail_bd2}), and noting that $X^{\alpha}_0$ is a $\{ 0,1\}$-valued random variable, we obtain that 
\begin{equation*}
\begin{aligned}
\E_0[X^{\alpha}_0] 
& = \E_0 \bigl[ X^{\alpha}_0 \bigm||\Phi_{(-\alpha R_L,0)}| \geq n_L \bigr] \prob( |\Phi_{(-\alpha R_L,0)}| \geq n_L  \\ & \qquad + \prob( |\Phi_{(-\alpha R_L,0)}|<n_L ) \E_0 \bigl[ X^{\alpha}_0 \bigm||\Phi_{(-\alpha R_L,0)}| < n_L \bigr] \\
& \leq \E_0 \bigl[ X^{\alpha}_0 \bigm||\Phi_{(-\alpha R_L,0)}| \geq n_L \bigr] + \prob( |\Phi_{(-\alpha R_L,0)}|<n_L )\\
& \leq L^{-\kappa}(1+o(1)).
\end{aligned}
\end{equation*}

Finally, the expected number of uncrossed gaps on $[0,L]$, which we denote $\E[\nucg]$, is obtained by integrating $\E_0[X_0]$ over $[0,L]$ with respect to the unit intensity measure of the Poisson process of node locations. Since $X^0\leq X^{\alpha}_0$, we conclude that 
$$
\E[\nucg]= L\E_0[X_0] \leq L^{1-\kappa}(1+o(1)).
$$
As $\kappa>0$ is arbitrary, we see by choosing $\kappa>1$ that $\E[\nucg]$ tends to zero as $L$ tends to infinity. Thus, the probability of an uncrossed gap is vanishing in the limit $L\to \infty$ for any choice of $\gamma>0$ and $R_L=\gamma \ln L$. This implies in particular that, at the critical scale of the connection function at which isolated nodes begin to appear, the probability of seeing an uncrossed gap is negligible. In other words, the primary mechanism responsible for causing disconnection in 1-D soft RGGs is the isolation of individual nodes.

As this is in sharp contrast to the situation for 1-D hard RGGs, in which disconnection is always due to uncrossed gaps, we comment briefly on the intuition behind the result. When the connection function is scaled as above, the isolation of a node is primarily determined by what happens in an interval of order $R_L$ around it. In the case of a hard RGG, this interval has to be empty. In the case of a soft RGG, it is not empty, but contains fewer nodes than expected, only $\epsilon R_L$, for some $\epsilon$ much smaller than 1. Isolation is then achieved by all $\epsilon R_L$ of these potential edges being absent. Both these events (having fewer nodes in the interval, and having no edges to these nodes) have probability which is exponentially small in $R_L$; at the scaling $R_L=\gamma \ln L$, this works out to a probability of order $1/L$ (of isolation per node). However, for there to be an uncrossed gap at the index node, $(\epsilon R_L/2)^2$ potential edges between the $\epsilon R_L/2$ nodes on each side of the index node have to be absent. This event has probability decaying exponentially in $R_L^2$, and is hence much less likely than the event of node isolation.

The intuition sketched out above is fleshed out in the appendix. It is shown in Appendix \ref{sec:app_conditional} that, conditional on a node $x$ being isolated, the probability that it also marks an uncrossed gap is vanishing in the above scaling regime. The scaling of the connection range at which uncrossed gaps appear is calculated (non-rigorously) in Appendix \ref{sec:app_scaling}, and shown to be much shorter than for isolated nodes. It is also shown there that, as the connection function tends towards that of the hard RGGs, the scaling regime for uncrossed gaps approaches that for isolated nodes.

%% file: Section8_Conclusion.tex
\section{Concluding Remarks}
\label{sec:conclusions}

In this paper, we proposed a model for soft random geometric graphs in one dimension in Section~\ref{sec:model}, and studied its connectivity properties. Our random graph model is restricted to point sets generated according to a homogenous Poisson process, but allows for very general connection functions.

In Section~\ref{sec:analysis}, we presented empirical evidence from simulations that the main modes of disconnection of our graphs are the presence of isolated nodes, and the presence of uncrossed gaps which partition the network into two or more disjoint clusters. Motivated by this evidence, in the next two sections we provided a rigorous mathematical analysis of the probabilities of isolated nodes or uncrossed gaps being present, in a suitable asymptotic scaling regime.

The analysis in Section~\ref{sec:iso} showed that the probability of the occurrence of isolated nodes shows an abrupt transition from nearly 1 to nearly 0 as the mean node degree increases. The transition is insensitive to the connection function, and occurs when the mean node degree is equal to the natural logarithm of the mean number of nodes.

The analysis in Section~\ref{sec:ucg} showed that uncrossed gaps have negligible probability in the scaling regime in which isolated nodes exhibit their transition. In other words, if any uncrossed gaps are present, then isolated nodes are present in abundance, but not conversely. Consequently, between the two, it is the presence of isolated nodes that is pivotal for determining connectivity. This was shown, not just for specific connection functions, but any function satisfying the mild assumptions of monotonocity, integrability, and unbounded support.

The finding that connectivity is determined by isolated nodes rather than uncrossed gaps is in stark contrast to hard random geometric graphs (the Gilbert or hard disk model), where it is uncrossed gaps that dominate. This does not contradict our analysis, since the Boolean connection function has bounded support. Nevertheless, it motivates the question of what happens with connection functions which approximate the Boolean threshold function, while having unbounded support. Simulation results presented in Section \ref{sec:analysis} (Figure \ref{WaxmanOverall}) show that isolated nodes are more prevalent for Waxman and uncrossed gaps for Rayleigh connection functions; the latter are closer to the Boolean model. Thus, for the Rayleigh connection function, and the system sizes simulated in this paper, which are quite large, it appears that our limiting analysis of uncrossed gaps does not apply.

This discrepancy between theory and simulations points to the need for an analysis of large but finite systems, to supplement the limit analysis presented in this paper. It is one of the prominent open problems raised by this work.

For the case of isolated nodes, we conjectured such a result in Section \ref{sec:iso}, in the form of a Poisson approximation for the number of isolated nodes. We presented simulation results in Figure \ref{fig:iso_wax} supporting this conjecture. A rigorous proof of this conjecture is an open problem, as is a similar result for uncrossed gaps. Finally, it also remains to be demonstrated rigorously that node isolation is the dominant mechanism for disconnection in the large system limit, not only in comparison to uncrossed gaps, but to all other possible forms of disconnection. Further open problems are to extend the analysis to more general point processes for modelling the node locations, and to temporal networks, namely those evolving dynamically over time, either through node mobility or through the formation and dissolution of links.

%% file: Section9_Ack.tex
\section*{Acknowledgements}

MW acknowledges the financial support from the EPSRC funded Centre for Doctoral Training in Communications (EP/I028153/1 and EP/L016656/1) and CD from the Spatially Embedded Networks project (EP/N002458/1). The authors had useful discussions with Lorenzo Federico and Alexander Kartun-Giles.

%% file: Section_Appendix.tex

\appendix

\section{Conditional probability of uncrossed gaps}
\label{sec:app_conditional}

We now look at the intuition behind why isolated nodes dominate uncrossed gaps as the cause of disconnection in soft RGGs, whereas the opposite is true of hard RGGs. We do this by conditioning on there being an isolated node at the origin, and computing the conditional probability that it also marks an uncrossed gap. In order to simplify calculations, we work with the Poisson process on the infinite real line.

Denote by $\prob_0$ the Palm probability corresponding to the presence of a point at the origin, and by $\hat{\prob}_0$ the measure obtained by conditioning further on this point being isolated. The corresponding expectations are denoted $\E_0$ and $\hat{\E}_0$. Then, under $\hat{\prob}_0$, the remaining points constitute an inhomogenous Poisson process, with intensity function 
$$
\mu(x) = 1-h(0,x) = 1-H(|x|), \; x\in \Reals.
$$
We now follow the same approach as in the previous section. Conditional on $\Phi_{(-\infty,0)}$, the point process restricted to the negative real line, the set of points on $(0,\infty)$ which are connected to some point in $\Phi_{(-\infty,0)}$ constitute an inhomogenous Poisson process with intensity 
\begin{equation}
    \label{eq:rate_cond_isol}
    \Lambda(x)=(1-H(x))\Bigl( 1-\exp\bigl( \sum_{z\in \Phi_{(-\infty,0)}} g(z,x) \bigr) \Bigr),
\end{equation}
where $g(x,y)=\ln(1-h(x,y))$. The total number of points is thus a Poisson random variable, with mean equal to $\int_0^{\infty} \Lambda(x) dx$. The event that the origin marks an uncrossed gap is the event that this Poisson random variable takes the value zero, which has probability 
\begin{equation}
    \label{eq:ucg0_cond_prob}
    \begin{aligned}
    \prob(\mbox{0 is uncrossed}) &= \hat{\E}_0 \Bigl[ \exp \Bigl( -\int_0^{\infty} \Lambda(x)dx \Bigr) \Bigr] \\
    &\geq \exp \Bigl( -\int_0^{\infty} \hat{\E}_0 [\Lambda(x)] dx \Bigr).
    \end{aligned}
\end{equation}
We have used Jensen's inequality to obtain the inequality above.

Now, by invoking Campbell's formula for the Laplace functional of the inhomogenous Poisson process $\Phi_{(-\infty,0)}$ \cite[Theorems 4.6 and 4.9]{haenggi2012stochastic}, we obtain from $\eqn$ (\ref{eq:rate_cond_isol}) that, for $x>0$, 
\begin{equation}
    \label{eq:cond_mean_rate}
    \begin{aligned}
    \hat{\E}_0[\Lambda(x)] &= (1-H(x)) \Bigl( 1 - \exp \Bigl( \int_{-\infty}^0 \bigl( e^{g(z,x)}-1) \bigr) (1-h(0,z))dz \Bigr) \Bigr) \\
    &= (1-H(x)) \Bigl( 1 - \exp \int_{-\infty}^0 -h(z,x)(1-h(0,z)) dz \Bigr) \\
    &= (1-H(x)) \Bigl( 1 - \exp \Bigl( -\int_0^{\infty} (1-H(z))H(x+z)dz \Bigr) \Bigr).
    \end{aligned}
\end{equation}
In the scaling regime in which the connection function $H(\cdot)$ is replaced by $H^L(\cdot)=H(\cdot/R_L)$, we can rewrite the above as 
$$
\hat{\E}[\Lambda(x)] = (1-H(x/R_L)) \Bigl( 1-\exp \Bigl( -R_L \int_{0}^{\infty} (1-H(z))H((x/R_L)+z)dz \Bigr) \Bigr),
$$
and so 
\begin{equation}
    \label{eq:connected_pts_mean_scaled}
    \int_0^{\infty} \hat{\E}[\Lambda(x)]dx = R_L \int_0^{\infty} (1-H(y)) \bigl( 1- e^{ -R_L \int_0^{\infty} (1- H(z))H(y+z)dz } \bigr)dy.
\end{equation}

Define $H^{-1}(x)=\inf \{ y\geq 0 : H(y) \leq x \}$, and note that integrable. By the assumption that $H$ is monotone decreasing, $1-H(z) \geq 1-x$ for all $z>H^{-1}(x)$. Hence,
$$
\int_0^{\infty} (1-H(z))H(y+z) dz \geq \frac{1}{2} \int_{H^{-1}(1/2)}^{\infty} H(y+z)dz.
$$
Substituting this in $\eqn$ (\ref{eq:connected_pts_mean_scaled}), invoking the inequality $1-H(y) \geq 1/2$ for all $y\geq H^{-1}(1/2)$ once more, and using the monotonicity of $H(\cdot)$, we get 
\begin{equation}
    \label{eq:connected_pts_mean_scaled_2}
    \begin{aligned}
    \int_0^{\infty} \hat{\E}[\Lambda(x)]dx  &\geq R_L \int_{H^{-1}(1/2)}^{\infty} \frac{1}{2} \bigl( 1- e^{ -(R_L/2) \int_{H^{-1}(1/2)}^{\infty} H(y+z)dz } \bigr)dy \\
    &\geq \frac{R_L}{2} \int_{H^{-1}(1/2)}^{2H^{-1}(1/2)} \Bigl( 1-\exp \Bigl(-\frac{R_L}{2} \int_{3H^{-1}(1/2)}^{\infty} H(z)dz \Bigr) \Bigr)dy.
    \end{aligned}
\end{equation}
Now, the integral in the exponent is a strictly positive constant, by the assumption that $H$ has unbounded support. Hence,
$$
\exp \Bigl(-\frac{R_L}{2} \int_{3H^{-1}(1/2)}^{\infty} H(z)dz \Bigr) \to 0 \mbox{ as } R_L\to \infty,
$$
and it follows that 
\begin{equation}
\begin{aligned}
    \label{eq:connected_pts_mean_scaled_3}
    \liminf_{R_L\to \infty} \frac{1}{R_L} \int_0^{\infty}\hat{\E}[\Lambda(x)]dx \geq \gamma &= \frac{1}{2} \int_{H^{-1}(1/2)}^{2H^{-1}(1/2)} 1 dy \\ 
    & = \frac{1}{2}H^{-1}(1/2).
\end{aligned}
\end{equation}

Since $H$ has unbounded support, $\gamma$ is a strictly positive constant, and we conclude that $\int_0^{\infty}\hat{\E}[\Lambda(x)]dx $, the expected number of nodes in $[0,\infty)$ which have a neighbour in $(-\infty,0)$, tends to infinity as $R_L$ tends to infinity. This does not prove that there is at least one such node with high probability, but it is at least strongly suggestive of it. Thus, $\eqn$ (\ref{eq:connected_pts_mean_scaled_3}) gives us strong reason to believe that the point at the origin, which was conditioned to be isolated, has very small probability of marking an uncrossed gap. This provides some partial intuition for why uncrossed gaps are rare in soft RGGs at the scaling threshold for the emergence of isolated nodes, whereas they are more prevalent that isolated nodes in hard RGGs.

\section{Scaling for emergence of uncrossed gaps}
\label{sec:app_scaling}

In Appendix \ref{sec:app_conditional}, we provided intuition for why uncrossed gaps are rare in the scaling regime at which isolated nodes emerge in 1-D soft RGGs. Here, we seek to identify the scaling regime at which uncrossed gaps emerge. We begin by calculating a bound on the expected number of uncrossed gaps within an interval $[0,L]$, in a soft RGG whose nodes are placed according to a unit rate PPP on the infinite real line; edges are then created independently, with probability $H(r)$ for nodes that are distance $r$ apart. We suppose that Assumption B from Section~\ref{sec:ucg} continues to hold.



Define $X_0$ to be the indicator that there is a point at the origin and that it marks an uncrossed gap, i.e., there are no edges between $\Phi_{(-\infty, 0)} \cup \lbrace 0 \rbrace$ and $\Phi_{(0, \infty)}$. We wish to calculate $\E_0[X_0]$, where $\mathbb{E}_0$ denotes expectation under the Palm measure conditional on the PPP having a point at the origin; thus, $\E_0[X_0]$ is the probability that the point at the origin constitutes an uncrossed gap. The expected number of uncrossed gaps in $[0,L]$ is then given by $\mathbb{E}[\nucg] = \int_0^L \mathbb{E}_0[X_0] dx = L\mathbb{E}_0[X_0]$.

The set of points on $(0,\infty)$ which have an edge to some point in $(-\infty,0]$ constitute a Cox process, which we denote by $N(\cdot)$; conditional on $\Phi_{(-\infty,0)}$, they constitute a Poisson process with intensity 
\begin{equation} \label{eq:Lambda_def}
\Lambda(y)= 1-(1-h(0,y))\prod_{z\in \Phi_{(-\infty,0)}}(1-h(z,y)) = 1-\exp \bigl( g(0,y)+\mbox{$\sum_{z\in \Phi_{(-\infty,0)}} g(z,y)$} \bigr),
\end{equation}
where $g(x,y)=\log(1-h(x,y))$. We use a capital letter to denote the intensity to make it explicit that the intensity is random, as it is a function of the random measure $\Phi$. 

Now, the event that there are no edges from $(-\infty,0]$ to $(0,\infty)$ is precisely the event that the Cox process of points in $(0,\infty)$ reached by such edges is empty, i.e., that $|N((0,\infty))|=0$. This is exactly the event whose indicator we defined as $X_0$. Hence, we obtain using the tower rule (law of iterated expectation) that 
\begin{equation} \label{eq:prob_split_0}
\E_0[X_0]= \E_0[ \E_0[X_0|\Phi_{(-\infty,0)}] ] = \E_0 \Bigl[\exp \Bigl( -\int_0^{\infty} \Lambda(y)dy \Bigr)\Bigr],
\end{equation}
where $\Lambda(\cdot)$ is given by $\eqn$ (\ref{eq:Lambda_def}).

A lower bound on this expectation is easy to compute using Jensen's inequality, which states that, if $X$ is a random variable and $f$ is a convex function, then
$$
f(\E [X]) \leq \E (f(X)).
$$
Applying this to $\eqn$ (\ref{eq:prob_split_0}) with $f(x)=e^{-x}$, we get
\begin{equation} \label{eqn:E_X_0}
    \E_0[X_0] \geq \exp \Bigl( -\E_0\Bigl[ \int_0^{\infty} \Lambda(y)dy \Bigr]\Bigr) = \exp \left(- \int_0^{\infty} \E[\Lambda (y)] dy  \right),
\end{equation}
where the interchange of integral and expectation is justified by Tonelli's theorem. We have also replaced $\E_0$ by $\E$, as the conditioning on having a point of the PPP at the origin has been taken into account in the expression for $\Lambda$.

The expectation of $\Lambda(y)$ can now be calculated using Campbell's formula for the characteristic functional of the Poisson point process \cite[Definition~4.7]{haenggi2012stochastic}. Letting $g(x,y)=\ln(1-h(x,y))$, we obtain from $\eqn$ (\ref{eq:Lambda_def}) that
\begin{equation}
    \begin{aligned}
    \E[\Lambda (y)] &= 1-(1-h(0,y))\E\bigl[ \exp\bigl( \mbox{$\sum_{z\in \Phi_{(-\infty,0)}} g(z,y)$} \bigr)\bigr] \\
    & = 1 - (1-h(0,y)) \exp \left(\int_{-\infty}^0\left( e^{g(z, y)} - 1 \right) dz\right) \\
    & = 1 - (1 - h(0, y))\exp\left(-\int_{-\infty}^0 h(z, y) dz \right) \\
    & = 1 - (1 - H(y)) \exp \left(-\int_y^{\infty} H(z) dz \right).
    \end{aligned}
\end{equation}
Integrating this over $y \in (0, \infty)$, we obtain
\begin{equation}\label{eqn:E_lambda_int}
    \begin{aligned}
    \int_0^{\infty} \E[\Lambda (y)] dy &= \int_0^{\infty} 1 - (1 - H(y)) \exp \left(-\int_y^{\infty} H(z) dz \right) dy \\
    & = \int_0^{\infty} \left(1 - e^{-\int_y^{\infty} H(z) dz} \right) dy + \int_0^{\infty} H(y) e^{-\int_y^{\infty} H(z) dz} dy \\
    & = \int_0^{\infty} \left(1 - e^{-\int_y^{\infty} H(z) dz} \right) dy + 1 - e^{-\int_0^{\infty} H(z) dz}.
    \end{aligned}
\end{equation}
Now, considering the scaled version of the connection function, $H^L(z)=H(z/R_L)$, we can rewrite the above as
\begin{equation}\label{eqn:E_lambda_int_scaled}
    \begin{aligned}
    \int_0^{\infty} \E[\Lambda (y)] dy 
    & = \int_0^{\infty} \left(1 - e^{-\int_y^{\infty} H^L(z) dz} \right) dy + 1 - e^{-\int_0^{\infty} H^L(z) dz} \\
    & = \int_0^{\infty} \left(1 - e^{-\int_y^{\infty} H(z/R_L) dz} \right) dy + 1 - e^{-\int_0^{\infty} H(z/R_L) dz} \\
    & \leq \int_0^{\infty} \left(1 - e^{-R_L\int_{y/R_L}^{\infty} H(z)dz} \right) dy + 1.
    \end{aligned}
\end{equation}

Note that $1-e^{-x}\leq x$ for all $x\in \Reals$, and so,
\begin{equation*} 
1 - e^{-R_L\int_{y/R_L}^{\infty} H(z)dz} \leq \min \Bigl\{ 1, R_L\int_{y/R_L}^{\infty} H(z)dz \Bigr\}.
\end{equation*}
Hence, for arbitrary $x>0$, we have
\begin{equation}
    \label{eq:integral_bd}
    \begin{aligned}
    \int_0^{\infty} 1 - e^{-R_L\int_{y/R_L}^{\infty} H(z)dz} dy 
    &\leq R_L x + \int_{R_L x}^{\infty} R_L\int_{y/R_L}^{\infty} H(z)dz dy \\
    &= R_L x + R_L \int_x^{\infty} \int_{R_L x}^{R_L z} H(z) dydz \\
    &\leq R_L x + R_L^2 \int_x^{\infty} zH(z)dz.
    \end{aligned}
\end{equation}

We shall henceforth restrict attention to the generalised Rayleigh connection function, $H(r)=\beta \exp(-(r/r_c)^{\eta})$, for a fixed $\eta>0$. 
For this connection function, we rewrite the last integral on the RHS of $\eqn$ (\ref{eq:integral_bd}) as 
\begin{equation}
    \label{eq:integral_tail}
    \begin{aligned}
    \int_x^{\infty} z H(z)dz &= \int_x^{\infty} \beta z e^{-(z/r_c)^{\eta}}dz \\ 
    &= \frac{\beta r_c^2}{\eta} \int_{\left(\frac{x}{r_c}\right)^{\eta}}^{\infty} u^{\frac{2}{\eta}-1} e^{-u} du \\
    &= \frac{\beta r_c^2}{\eta} \Gamma \left( \frac{2}{\eta}, \left(\frac{x}{r_c}\right)^{\eta} \right),
    \end{aligned}
\end{equation}
where the second line comes from the substitution $u = \left( \frac{z}{r_c} \right)^{\eta}$ and $\Gamma(\alpha,y)$ denotes the incomplete Gamma function,
$$
\Gamma(\alpha,y) = \int_y^{\infty} z^{\alpha-1}e^{-z}dz.
$$
Observe that
\begin{align*}
\frac{\Gamma(\alpha,y)}{y^{\alpha-1}e^{-y}} &= \int_y^{\infty} \Bigl( \frac{x}{y} \Bigr)^{\alpha-1} e^{-(x-y)}dx \\
&= \int_0^{\infty} \Bigl( 1+\frac{z}{y} \Bigr)^{\alpha-1} e^{-z} dz,
\end{align*}
which tends to 1 as $y$ tends to infinity, by the Dominated Convergence Theorem. Hence, taking $x=r_c(\ln R_L)^{1/\eta}$ in $\eqn$ (\ref{eq:integral_tail}), we see that 
\begin{equation*}
    \int_{r_c (\ln R_L)^{1/\eta}}^{\infty} zH(z)dz = \frac{\beta r_c^2}{\eta} \Gamma\left(\frac{2}{\eta}, \ln R_L \right)\sim \frac{\beta r_c^2}{\eta} \frac{(\ln R_L)^{\frac{2}{\eta}-1}}{R_L} \mbox{ as $R_L\to \infty$}.
\end{equation*}
Here, for functions $f$ and $g$ on $\Reals_+$, we write $f(x)\sim g(x)$ as $x\to \infty$ to denote that $f(x)/g(x)$ tends to 1 as $x$ tends to infinity. Substituting the above in $\eqn$ (\ref{eq:integral_bd}), we get 
\begin{equation*}
    \int_0^{\infty} 1 - e^{-R_L\int_{y/R_L}^{\infty} H(z)dz} dy \leq r_c R_L (\ln R_L)^{1/\eta} + (1+o(1)) \frac{\beta r_c^2}{\eta} R_L (\ln R_L)^{\frac{2}{\eta}-1}.
\end{equation*}

Combining the above with $\eqns$ (\ref{eqn:E_X_0}) and (\ref{eqn:E_lambda_int_scaled}), we obtain the following lower bound on the probability that there is an uncrossed gap at the origin:
\begin{equation*}
    \E_0[X_0] \geq \exp \Bigl( -C R_L (\ln R_L)^{\max \{ \frac{1}{\eta}, \frac{2}{\eta}-1 \} } \Bigr),
\end{equation*}
where $C>0$ is a fixed constant that does not grow with $R_L$. 

The expected number of uncrossed gaps in $[0,L]$ is given by $\E[\nucg]=L\E_0[X_0]$. Hence, it follows from the equation above that 
\begin{equation} \label{eq:nucg_lower_bd_2}
\E[\nucg] \geq L\exp \bigl( -C R_L (\ln R_L)^{\theta} \bigr) \mbox{ where } \theta = \max \left\{ \frac{1}{\eta}, \frac{2}{\eta}-1 \right\}.
\end{equation}
The above expression motivates us to consider the scaling regime 
\begin{equation}
    \label{eq:ucg_crit_scaling}
    R_L = \gamma \frac{\ln L}{(\ln \ln L)^{\theta}}.
\end{equation}
A straightforward calculation shows that 
\begin{align*}
CR_L (\ln R_L)^{\theta} &= \gamma C \ln L \Bigl( 1+ \frac{\theta \ln \ln \ln L + \ln \gamma}{\ln \ln L} \Bigr)^{\theta} \\
&= \gamma C \ln L (1+o(1)).
\end{align*}
If we conjecture that the inequality in $\eqn$ (\ref{eq:nucg_lower_bd_2}) is an approximate equality, then it follows that the expected number of uncrossed gaps exhibits a sharp threshold at $\gamma=1/C$, in the sense that
\begin{equation*}
\E[\nucg] \to \begin{cases}
0, & \gamma>1/C, \\
+\infty, & \gamma<1/C, 
\end{cases}
\end{equation*}
as $L$ tends to infinity. 

Thus, our calculations lead us to conjecture that uncrossed gaps appear when the connection range scales as $R_L = \ln L/(C(\ln \ln L)^{\theta})$; here $\theta$ is the related to the power law in the exponent, $\eta$ of the generalised Rayleigh connection function, while $C$ is a constant that depends in a complicated way on the parameters of the connection function. We now make a few remarks about this scaling regime. Firstly, $C$ and $\theta$ depend on details of the connection function, and are not universal, in contrast to the threshold for isolated nodes. Secondly, $R_L$ is significantly smaller than the $\ln L$ threshold for isolated nodes, confirming that uncrossed gaps appear only when connection functions are of much shorter range than required for the appearance of isolated nodes. Finally, the hard RGG model is a limiting case of the generalised Rayleigh connection function as $\eta$ tends to infinity; correspondingly, $\theta=1/\eta$ tends to zero, and the threshold value of $R_L$ tends to a constant multiple of $\ln L$. Thus, we recover the $\ln L$ scaling for the emergence of uncrossed gaps in the hard RGG in the limit.

We also comment here on the statements made at the end of Section \ref{sec:analysis} regarding the fact that we don't expect to see these results in simulations unless very large system sizes are analysed. For clarity, we will denote by $R^{ucg}_L$ the scaling regime in which there is a sharp threshold for the appearance of uncrossed gaps, and we denote by $R^{iso}_L$ the relevant scaling regime for isolated nodes. These have been shown in this paper to be

$$
R_L^{ucg} = \frac{\ln L}{C (\ln \ln L)^{\theta}},
$$

and

$$
R_L^{iso} = \frac{\ln L}{2 \|H\|_1}.
$$

Therefore, to find the value of $L$ (i.e the system size) at which there is a crossover between uncrossed gaps being the more important factor when discussing connectivity and when isolated nodes are more important, we need to find the value of $L$ for which $R_L^{ucg} = R_L^{iso}$, which we denote by $L^{*}$. This occurs at 

$$
\frac{\ln L^{*}}{2 \|H\|_1} = \frac{\ln L^{*}}{C (\ln \ln L^{*})^{\theta}}.
$$

Rearranging this we see that 

$$
L^{*} = \exp\left(\exp\left(\left(\frac{2 \|H\|_1}{C}\right)^{1 / \theta}\right)\right).
$$

For $\eta \geq 1$, $\theta = \frac{1}{\eta}$ and hence this becomes

\begin{equation}
\label{eqn:critical_L}
    L^{*} = \exp\left(\exp\left(\left(\frac{2 \|H\|_1}{C}\right)^{\eta}\right)\right).
\end{equation}

The calculation of the value of $C$ is an interesting problem in its own right and is hence left as future work. However, Eqn. (\ref{eqn:critical_L}) is able to tell us that $L^{*}$ is delicately related to this value of $C$, and could be extremely large.